# DYNAMICAL LARGE DEVIATIONS FOR THE BOUNDARY DRIVEN WEAKLY ASYMMETRIC EXCLUSION PROCESS

By Lorenzo Bertini, Claudio Landim and Mustapha Mourragui

*Università di Roma "La Sapienza," IMPA and Université de Rouen and Université de Rouen*

We consider the weakly asymmetric exclusion process on a bounded interval with particle reservoirs at the endpoints. The hydrodynamic limit for the empirical density, obtained in the diffusive scaling, is given by the viscous Burgers equation with Dirichlet boundary conditions. We prove the associated dynamical large deviations principle.

**1. Introduction.** The analysis of the large deviations is asymptotic as the number of degrees of freedom diverges, for the *stationary nonequilibrium states* of interacting particle systems have recently proved to be an important step in the physical description of such states and a rich source of mathematical problems. Referring to [4, 8] for two recent reviews on this topic, we briefly outline its basic points. We discuss only *stochastic lattice gases* for which the underlying random fluctuations ensure the necessary ergodicity for a rigorous analysis. The stationary nonequilibrium states are characterized by a flow of mass through the system and the corresponding dynamics are not reversible. The main difference with respect to (reversible) equilibrium states follows. In equilibrium the invariant measure, which determines the thermodynamic properties, is given for free by the Gibbs distribution specified by the Hamiltonian. On the contrary, in nonequilibrium states the construction of the stationary state requires the solution of a dynamical problem.

Since we are interested only in the macroscopic description, only the thermodyamic observables are relevant. For lattice gases there is only one of such observable, which is the *empirical density*. In equilibrium states, the thermodynamic functional, like the free energy, can then be identified [6, 22, 24] with the large deviation rate function for the empirical density when the particles are distributed according to the Gibbs measure. Provided one replaces the Gibbs measure with the invariant measure (which is in general not explicitly known), the above statement is meaningful also for stationary nonequilibrium states, and it is the definition of nonequilibrium free energy adopted in [4, 8].

A typical generic feature of nonequilibrium stationary states is the presence of long-range correlations. The large deviation rate functional, which has been identified with the thermodynamic functional, is nonlocal. In this respect nonequilibrium stationary states behave quite differently from equilibrium states. As it has been shown in concrete examples, the nonequilibrium rate functional is also not necessarily convex.

In absence of the Gibbs principle, a basic problem is the characterization of the nonequilibrium free energy in concrete simple models of stochastic lattice gases. While the approach discussed in [8] is based on powerful combinatorial methods which exploit the special feature of the exclusion processes, we recall the dynamical/variational approach reviewed in [4]. One first fixes a macroscopic time interval $[0,T]$ and analyzes the dynamical behavior of the empirical density over such an interval. The law of large numbers for the empirical density is then called *hydrodynamic limit*, and, in the context of diffusive scaling limit here considered, it is given by a parabolic evolution equation. The next step is the proof of the associated dynamical large deviation principle, namely, to compute the asymptotic probability of observing a given large fluctuation in the dynamics of the empirical density. This asymptotic probability can be expressed by a suitable rate functional on the set of space-time trajectories. Finally, one minimizes the dynamical rate functional on all the paths starting from the stationary profile, that is, the stationary solution of the hydrodynamic equation ending on a fixed profile. The solution of this variational problem then coincides with the nonequilibrium free energy.

The exclusion process is a very simple lattice gas: the only interaction is due to the exclusion condition. A particle can therefore jump to its neighboring sites, but the jump takes place only if the arrival site is not occupied. We consider this process on the bounded lattice $[-N+1, N-1] \cap \mathbb{Z}$, $N \geq 1$, in contact with particles' reservoirs at the endpoints, so that to the bulk dynamics we add birth and death processes at the sites $\pm(N-1)$. In the case of the boundary driven symmetric exclusion process characterized by symmetric bulk jump rates, the program outlined has been rigorously implemented in [3]. Of course, the nonequilibrium free energy functional thus obtained coincides with the one deduced in [9] by combinatorial methods.

Here we analyze instead boundary driven weakly asymmetric exclusion processes in which the asymmetry of the bulk rates is of order $E/N$ for some fixed $E \in \mathbb{R}$. For this model on the whole lattice, the hydrodynamic limit has been proven in [7, 16, 20], and the hydrodynamic equation is the viscous Burgers equation (see [14, 15, 18] for the hydrodynamic limit of boundary driven models). Referring to [5, 13] for the computation of the nonequilibrium free energy functional, in this paper we prove the dynamical large deviation principle associated to the hydrodynamic limit. The general methods developed in [20] for the symmetric exclusion process, and adapted in [3] to the boundary driven case, are applied with simple modifications, but there is a somewhat delicate technical point. The basic strategy in the proof of the lower bound consists in obtaining this bound for smooth paths and then applying the following density argument. Given a path $\pi$ with a finite rate functional, that is, $I(\pi) < \infty$, one constructs a suitable sequence of smooth paths $\{\pi_n\}$ such that $\pi_n \to \pi$ and $I(\pi_n) \to I(\pi)$. By the lower semicontinuity of $I$ we have $\liminf_n I(\pi_n) \geq I(\pi)$ for *any* sequence $\{\pi_n\}$, but one needs to show that the equality actually holds for a *suitable* sequence $\{\pi_n\}$. The proof of this step in the symmetric case in [3, 20] takes advantage of the convexity of the functional $I$; this property does not, however, hold for the weakly asymmetric exclusion process. Following [26, 27], we modify the definition of the rate functional $I$ requiring that a suitable energy estimate holds. As shown here in detail, this energy estimate provides the necessary compactness to carry out the above density argument. The arguments needed to prove this density are essentially an adaptation to the boundary driven case of those developed in [26, 27]. The lack of translation invariance requires new tools.

The modification in the definition of the rate functional $I$ makes the proof of the upper bound harder than the one in [3, 20]: one needs to show that the energy estimate holds with probability superexponentially close to one. This step is also discussed here in detail. A similar approach has been followed in [23] for particle systems with Kac interaction and random potential.

The proof of a dynamical large deviations principle for the empirical measure of boundary driven interacting particle systems is the first step in the derivation of the nonequilibrium free energy, a thermodynamical functional of considerable interest in mathematical physics. Based on the results presented here, we obtain in [5] the nonequilibrium free energy of weakly asymmetric exclusion processes and show that its limit, as the asymmetry diverges, $\Gamma$-converges to the nonequilibrium free energy of the asymmetric exclusion process obtained in [10] through combinatorial methods.

## 2. Notation and results.

*The boundary driven weakly asymmetric exclusion process.* Fix an integer $N \geq 1$, $E \in \mathbb{R}$, $0 < \rho_- \leq \rho_+ < 1$ and let $\Lambda_N := \{-N+1, \ldots, N-1\}$. The configuration space is $\Sigma_N := \{0,1\}^{\Lambda_N}$; elements of $\Sigma_N$ are denoted by $\eta$ so that $\eta(x) = 1$, resp. 0, if site $x$ is occupied, resp. empty, for the configuration $\eta$. We denote by $\sigma^{x,y}\eta$ the configuration obtained from $\eta$ by exchanging the occupation variables $\eta(x)$ and $\eta(y)$, that is,

$$(\sigma^{x,y}\eta)(z) := \begin{cases} \eta(y), & \text{if } z = x, \\ \eta(x), & \text{if } z = y, \\ \eta(z), & \text{if } z \neq x, y, \end{cases}$$

and by $\sigma^x\eta$ the configuration obtained from $\eta$ by flipping the configuration at $x$, that is,

$$(\sigma^x\eta)(z) := \begin{cases} 1 - \eta(x), & \text{if } z = x, \\ \eta(z), & \text{if } z \neq x. \end{cases}$$

The one-dimensional boundary driven weakly asymmetric exclusion process is the Markov process on $\Sigma_N$ whose generator $L_N$ can be decomposed as

$$L_N = L_{0,N} + L_{-,N} + L_{+,N},$$

where the generators $L_{0,N}$, $L_{-,N}$, $L_{+,N}$ act on functions $f : \Sigma_N \to \mathbb{R}$ as

$$(L_{0,N}f)(\eta) = \frac{N^2}{2} \sum_{x=-N+1}^{N-2} e^{-E/(2N)[\eta(x+1)-\eta(x)]} [f(\sigma^{x,x+1}\eta) - f(\eta)],$$

$$(L_{-,N}f)(\eta) = \frac{N^2}{2} c_-(\eta(-N+1))[f(\sigma^{-N+1}\eta) - f(\eta)]$$

$$(L_{+,N}f)(\eta) = \frac{N^2}{2} c_+(\eta(N-1))[f(\sigma^{N-1}\eta) - f(\eta)],$$

where $c_\pm : \{0,1\} \to \mathbb{R}$ are given by

$$c_\pm(\zeta) := \rho_\pm e^{\mp E/(2N)}(1-\zeta) + (1-\rho_\pm)e^{\pm E/(2N)}\zeta.$$

Notice that the (weak) external field is $E/(2N)$, and, in view of the diffusive scaling limit, the generator has been speeded up by $N^2$. We denote by $\eta_t$ the Markov process on $\Sigma_N$ with generator $L_N$ and by $\mathbb{P}_\eta^N$ its distribution if the initial configuration is $\eta$. Note that $\mathbb{P}_\eta^N$ is a probability measure on the path space $D(\mathbb{R}_+, \Sigma_N)$, which we consider endowed with the Skorohod topology and the corresponding Borel $\sigma$-algebra. Expectation with respect to $\mathbb{P}_\eta^N$ is denoted by $\mathbb{E}_\eta^N$.

Since the Markov process $\eta_t$ is irreducible, for each $N \geq 1$, $E \in \mathbb{R}$, and $0 < \rho_- \leq \rho_+ < 1$ there exists a unique invariant measure $\mu_E^N$ in which we drop the dependence on $\rho_\pm$ from the notation. Let $\varphi_\pm := \log[\rho_\pm/(1-\rho_\pm)]$ be the chemical potential of the boundary reservoirs, and set $E_0 := (\varphi_+ - \varphi_-)/2$. A simple computation shows that if $E = E_0$ then the process $\eta_t$ is reversible with respect to the product measure

$$\mu_{E_0}^N(\eta) = \prod_{x=-N+1}^{N-1} \frac{e^{\overline{\varphi}_{E_0}^N(x)\eta(x)}}{1 + e^{\overline{\varphi}_{E_0}^N(x)}},$$

where

$$\overline{\varphi}_{E_0}^N(x) := \varphi_- \frac{N-x}{2N} + \varphi_+ \frac{N+x}{2N}.$$

On the other hand, for $E \neq E_0$ the invariant measure $\mu_E^N$ cannot be written in a simple form.

*The dynamical large deviation principle.* We denote by $u \in [-1,1]$ the macroscopic space coordinate and by $\langle \cdot, \cdot \rangle$ the inner product in $L_2([-1,1], du)$. We set

$$\mathcal{M} := \{\rho \in L_\infty([-1,1], du) : 0 \leq \rho \leq 1\},$$

which we equip with the topology induced by the weak convergence of measures, namely a sequence $\{\rho^n\} \subset \mathcal{M}$ converges to $\rho$ in $\mathcal{M}$ if and only if $\langle \rho^n, G \rangle \to \langle \rho, G \rangle$ for any continuous function $G : [-1,1] \to \mathbb{R}$. Note that $\mathcal{M}$ is a compact Polish space that we consider endowed with the corresponding Borel $\sigma$-algebra. The empirical density of the configuration $\eta \in \Sigma_N$ is defined as $\pi^N(\eta)$ where the map $\pi^N : \Sigma_N \to \mathcal{M}$ is given by

$$\pi^N(\eta)(u) := \sum_{x=-N+1}^{N-1} \eta(x) \mathbf{1}\left\{\left[\frac{x}{N} - \frac{1}{2N}, \frac{x}{N} + \frac{1}{2N}\right)\right\}(u),$$

in which $\mathbf{1}\{A\}$ stands for the indicator function of the set $A$. Let $\{\eta^N\}$ be a sequence of configurations with $\eta^N \in \Sigma_N$. If the sequence $\{\pi^N(\eta^N)\} \subset \mathcal{M}$ converges to $\rho$ in $\mathcal{M}$ as $N \to \infty$, we say that $\{\eta^N\}$ is *associated* with the macroscopic density profile $\rho \in \mathcal{M}$.

Given $T > 0$, we denote by $D([0,T]; \mathcal{M})$ the Skorohod space of paths from $[0,T]$ to $\mathcal{M}$ equipped with its Borel $\sigma$-algebra. Elements of $D([0,T], \mathcal{M})$ will be denoted by $\pi \equiv \pi_t(u)$ and sometimes by $\pi(t, u)$. Note that the evaluation map $D([0,T]; \mathcal{M}) \ni \pi \mapsto \pi_t \in \mathcal{M}$ is not continuous for $t \in (0,T)$ but is continuous for $t = 0, T$. We denote by $\pi^N$, also, the map from $D([0,T]; \Sigma_N)$ to $D([0,T]; \mathcal{M})$ defined by $\pi^N(\eta_\cdot)_t := \pi^N(\eta_t)$. The notation $\pi^N(t, u)$ is also used.

Fix a profile $\gamma \in \mathcal{M}$ and consider a sequence $\{\eta^N : N \geq 1\}$ associated to $\gamma$. Let $\eta_t^N$ be the boundary driven weakly asymmetric exclusion process starting from $\eta^N$. In [7, 16, 20] it is proven that as $N \to \infty$ the sequence of random variables $\{\pi^N(\eta_\cdot^N)\}$, which take values in $D([0,T], \mathcal{M})$ and converge in probability to the path $\rho \equiv \rho_t(u)$, $(t, u) \in [0, T] \times [-1, 1]$ which solves the viscous Burgers equation with Dirichlet boundary conditions at $\pm 1$, that is,

(2.1)
$$\begin{cases} \partial_t \rho + \frac{E}{2} \nabla \chi(\rho) = \frac{1}{2} \Delta \rho, \\ \rho_t(\pm 1) = \rho_\pm, \\ \rho_0(u) = \gamma(u), \end{cases}$$

where $\chi : [0,1] \to \mathbb{R}_+$ is the mobility of the system, $\chi(a) = a(1-a)$, and $\nabla$, resp. $\Delta$, denotes the derivative, resp. the second derivative, with respect to $u$. In fact the proof presented in [7, 16] is in real line, while the one in [20] is on the torus. The arguments, however, can be adapted to the boundary driven case (see [14, 15, 18] for the hydrodynamic limit of different boundary driven models).

The main result of this paper is the large deviations principle associated with the above law of large numbers. In order to state this result some more notation is required. For $T > 0$ and positive integers $m, n$, we denote by $C^{m,n}([0,T] \times [-1,1])$ the space of functions $G \equiv G_t(u) : [0,T] \times [-1,1] \to \mathbb{R}$ with $m$ derivatives in time, $n$ derivatives in space which are continuous up to the boundary. We improperly denote by $C_0^{m,n}([0,T] \times [-1,1])$ the subset of $C^{m,n}([0,T] \times [-1,1])$ of the functions which vanish at the endpoints of $[-1,1]$, that is, $G \in C^{m,n}([0,T] \times [-1,1])$ belongs to $C_0^{m,n}([0,T] \times [-1,1])$ if and only if $G_t(\pm 1) = 0$, $t \in [0,T]$.

Let the energy $\mathcal{Q} : D([0,T], \mathcal{M}) \to [0, \infty]$ be given by

$$\mathcal{Q}(\pi) = \sup_G \left\{ \int_0^T dt \int_{-1}^1 du\, \pi(t,u)(\nabla G)(t,u) - \frac{1}{2} \int_0^T dt \int_{-1}^1 du\, G(t,u)^2 \chi(\pi(t,u)) \right\},$$

where the supremum is carried over all smooth functions $G : [0,T] \times (-1,1) \to \mathbb{R}$ with compact support. In Section 4 we show that the energy $\mathcal{Q}$ is convex and lower semicontinuous. Moreover, if $\mathcal{Q}(\pi)$ is finite, $\pi$ has a generalized space derivative, $\nabla \pi$, and

$$\mathcal{Q}(\pi) = \frac{1}{2} \int_0^T dt \int_{-1}^1 du\, \frac{(\nabla \pi_t)^2}{\chi(\pi_t)}.$$

Fix a function $\gamma \in \mathcal{M}$ which corresponds to the initial profile. For each $H$ in $C_0^{1,2}([0,T] \times [-1,1])$, let $\hat{J}_H = \hat{J}_{T,H,\gamma} : D([0,T], \mathcal{M}) \to \mathbb{R}$ be the functional given by

$$\begin{aligned}
\hat{J}_H(\pi) := & \langle \pi_T, H_T \rangle - \langle \gamma, H_0 \rangle - \int_0^T dt \langle \pi_t, \partial_t H_t \rangle \\
& - \frac{1}{2} \int_0^T dt \langle \pi_t, \Delta H_t \rangle + \frac{\rho_+}{2} \int_0^T dt\, \nabla H_t(1) \\
& - \frac{\rho_-}{2} \int_0^T dt\, \nabla H_t(-1) \\
& - \frac{E}{2} \int_0^T dt \langle \chi(\pi_t), \nabla H_t \rangle - \frac{1}{2} \int_0^T dt \langle \chi(\pi_t), (\nabla H_t)^2 \rangle.
\end{aligned}$$

(2.2)

Let $\hat{I}_T(\cdot|\gamma) \colon D([0,T],\mathcal{M}) \longrightarrow [0,+\infty]$ be the functional defined by

$$\hat{I}_T(\pi|\gamma) := \sup_{H \in C_0^{1,2}([0,T] \times [-1,1])} \hat{J}_H(\pi). \tag{2.3}$$

The rate functional $I_T(\cdot|\gamma) \colon D([0,T],\mathcal{M}) \to [0,\infty]$ is given by

$$I_T(\pi|\gamma) = \begin{cases} \hat{I}_T(\pi|\gamma), & \text{if } \mathcal{Q}(\pi) < \infty, \\ \infty, & \text{otherwise.} \end{cases} \tag{2.4}$$

We prove in Theorem 4.2 that the functional $I_T(\cdot|\gamma)$ is lower semicontinuous and has compact level sets, and in Lemma 4.3 that any path $\pi$ with with finite rate function, $I_T(\pi|\gamma) < \infty$, is continuous in time and satisfies the boundary conditions $\pi(0,\cdot) = \gamma(\cdot)$, $\pi(\cdot,\pm 1) = \rho_\pm$. In Section 5 we show that any trajectory $\pi$ with finite rate function can be approximated by a sequence of smooth trajectories $\{\pi^n : n \geq 1\}$ such that $I_T(\pi^n|\gamma)$ converges to $I_T(\pi|\gamma)$. These properties of the rate function $I_T(\cdot|\gamma)$ hold in a general context described in Section 4.

The main result of this article reads as follows.

THEOREM 2.1. *Fix $T > 0$ and an initial profile $\gamma$ in $\mathcal{M}$. Consider a sequence $\{\eta^N : N \geq 1\}$ of configurations associated to $\gamma$. Then, the sequence of probability measures $\{\mathbb{P}_{\eta^N}^N \circ (\pi^N)^{-1} : N \geq 1\}$ on $D([0,T],\mathcal{M})$ satisfies a large deviation principle with speed $N$ and good rate function $I_T(\cdot|\gamma)$. Namely, $I_T(\cdot|\gamma) \colon D([0,T];\mathcal{M}) \to [0,\infty]$ has compact level sets and for each closed set $\mathcal{C} \subset D([0,T],\mathcal{M})$ and each open set $\mathcal{O} \subset D([0,T],\mathcal{M})$,*

$$\varlimsup_{N \to \infty} \frac{1}{N} \log \mathbb{P}_{\eta^N}^N(\pi^N \in \mathcal{C}) \leq -\inf_{\pi \in \mathcal{C}} I_T(\pi|\gamma),$$

$$\varliminf_{N \to \infty} \frac{1}{N} \log \mathbb{P}_{\eta^N}^N(\pi^N \in \mathcal{O}) \geq -\inf_{\pi \in \mathcal{O}} I_T(\pi|\gamma).$$

We provide in Section 4 an explicit formula for the rate function $I_T(\cdot|\gamma)$. Denote by $C_K^\infty((0,T) \times (-1,1))$ the infinitely differentiable functions $H \colon (0,T) \times (-1,1) \to \mathbb{R}$ with compact support. For a trajectory $\pi$ in $D([0,T],\mathcal{M})$, let $\mathcal{H}_0^1(\chi(\pi))$ be the Hilbert space induced by $C_K^\infty((0,T) \times (-1,1))$ endowed with the scalar product defined by

$$\langle\!\langle G, H \rangle\!\rangle_{1,\chi(\pi)} = \int_0^T dt \int_{-1}^1 du (\nabla G)(t,u)(\nabla H)(t,u)\chi(\pi(t,u)).$$

Induced means that we first declare two functions $F$, $G$ in $C_K^\infty((0,T) \times (-1,1))$ to be equivalent if $\langle\!\langle F - G, F - G \rangle\!\rangle_{1,\chi(\pi)} = 0$ and then we complete the quotient space with respect to the norm induced by the scalar product. Denote by $\|\cdot\|_{1,\chi(\pi)}$ the norm associated to the scalar product $\langle\!\langle \cdot, \cdot \rangle\!\rangle_{1,\chi(\pi)}$.

Let $\mathcal{H}^{-1}(\chi(\pi))$ be the dual of $\mathcal{H}_0^1(\chi(\pi))$. It is a Hilbert space equipped with the norm $\|\cdot\|_{-1,\chi(\pi)}$ defined by

$$\|L\|_{-1,\chi(\pi)}^2 = \sup_{G \in C_K^\infty((0,T)\times(-1,1))} \{2\langle\!\langle L, G\rangle\!\rangle - \|G\|_{1,\chi(\pi)}^2\}.$$

In this formula, $\langle\!\langle L, G\rangle\!\rangle$ stands for the value of the linear form $L$ at $G$. We prove in Section 4 that if $\pi$ is a trajectory with finite rate function, then

$$I_T(\pi|\gamma) = \tfrac{1}{2}\|\partial_t\pi - (1/2)\Delta\pi + (E/2)\nabla\chi(\pi)\|_{-1,\chi(\pi)}^2.$$

A large deviations principle for the symmetric simple exclusion process, $E = 0$, with periodic boundary conditions has been proved in [20]. It has been extended in [3] to symmetric exclusion processes in contact with reservoirs. In both cases the rate function is $\hat{I}_T(\cdot|\gamma)$.

By Lemma 2.1.1 in [11], there is uniqueness of rate functions in Polish spaces. In particular, for the symmetric simple exclusion process, $\hat{I}_T = I_T$. Equivalently, any path $\pi$ with finite rate function, $\hat{I}_T(\pi|\gamma) < \infty$, has finite energy.

**3. The large deviations principle.** We prove in this section, relying on some properties of the rate function that we prove later, the large deviations principle stated in Theorem 2.1.

The approach differs slightly from the original one in [12, 20] due to definition (2.4) of the rate function $I_T(\cdot|\gamma)$ which is set to be $+\infty$ on the set of paths $\pi$ with infinite energy $[\mathcal{Q}(\pi) = +\infty]$ [26, 27]. The rate function $I_T(\cdot|\gamma)$, being larger than $\hat{I}_T(\cdot|\gamma)$, the original one in [12, 20], the proof of the upper bound becomes harder and the one of the lower bound easier. As discussed in the introduction, this modification is needed for the following reason. In the lower bound part, one first proves the estimate for suitable "nice" trajectories and then one shows that any path with finite rate function can be approximated by a sequence of "nice" trajectories with convergence of the associated large deviations probability. The procedure used in [20] for this step relies strongly on the convexity of the rate functional which allows the approximation of a path by taking convolutions with a smooth ad-hoc kernel. However, the convexity of the rate function is a special feature of the symmetric exclusion process. Without such convexity, one is only able to approximate trajectories with finite energy, thus forcing the above redefinition of the rate function. The boundary conditions introduce a second obstacle which prevent convolutions with space-independent kernels since the rate function equals $+\infty$ for trajectories which do not meet the boundary conditions.

Denote by $\bar{\rho}$ the stationary density profile, that is, the unique solution of the elliptic equation

$$\begin{cases} E\nabla[\rho(1-\rho)] = \Delta\rho, \\ \rho(\pm 1) = \rho_\pm. \end{cases}$$

Denote by $\nu_N$ the product measure with density profile $\bar{\rho}$. The marginals of $\nu_N$ are given by

$$\nu_N\{\eta : \eta(x) = 1\} = \bar{\rho}(x/N), \qquad -N+1 \leq x \leq N-1.$$

3.1. *Superexponential estimates.* It is well known that one of the main steps in the derivation of a large deviation principle for the empirical density is a super-exponential estimate which allows the replacement of local functions by functionals of the empirical density in the large deviations regime. The problem consists of estimating expressions such as $\langle V, f^2 \rangle_{\mu^N}$ in terms of the Dirichlet form $\langle -L_N f, f \rangle_{\mu^N}$ where $V$ is a local function and $\langle \cdot, \cdot \rangle_{\mu^N}$ represents the inner product with respect to some probability measure $\mu^N$.

In the context of boundary-driven processes, the fact that the invariant measure is not known explicitly introduces a technical difficulty. Following [3, 21] we fix $\nu_N$, the product measure with density profile $\bar{\rho}$, as reference measure and estimate everything with respect to $\nu_N$. Note, however, that since $\nu_N$ is not the invariant measure, there are no reasons for $\langle -L_N f, f \rangle_{\nu_N}$ to be positive. The first statement shows that this expression is almost positive.

For a function $f : \Sigma_N \to \mathbb{R}$, let

$$D_{0,N}(f) = \sum_{x=-N+1}^{N-2} \int [f(\sigma^{x,x+1}\eta) - f(\eta)]^2 \, d\nu_N,$$

$$D_{\pm,N}(f) = \int [f(\sigma^{\pm(N-1)}\eta) - f(\eta)]^2 \, d\nu_N.$$

LEMMA 3.1. *There exist constants $A_0, C_0 > 0$ depending only on $\rho_\pm$, $E$, such that*

$$\langle L_{0,N} f, f \rangle_{\nu_N} \leq -A_0 N^2 D_{0,N}(f) + C_0 N \langle f, f \rangle_{\nu_N},$$
$$\langle L_{\pm,N} f, f \rangle_{\nu_N} \leq -A_0 N^2 D_{\pm,N}(f) + C_0 \langle f, f \rangle_{\nu_N}$$

*for all functions $f : \Sigma_N \to \mathbb{R}$.*

The proof of this lemma is elementary and left to the reader. The fact that $\bar{\rho}$ is the stationary profile is irrelevant. We may replace in the statement of the lemma the product measure $\nu_N$, associated with $\bar{\rho}$, by the product measure associated with any smooth profile with the correct boundary conditions at $\pm 1$. Lemma 3.1, together with the computation presented in [2], page 78, for nonreversible processes, allows one to prove the super-exponential estimates stated below in Theorem 3.2.

Given a cylinder function $\Psi$, that is a function on $\{0,1\}^{\mathbb{Z}}$ depending on $\eta(x)$, $x \in \mathbb{Z}$, only through finitely many $x$, denote by $\widetilde{\Psi}(\alpha)$ the expectation of $\Psi$ with respect to $\nu_\alpha$, the Bernoulli product measure with density $\alpha$:

$$\widetilde{\Psi}(\alpha) = E_{\nu_\alpha}[\Psi].$$

Denote by $\{\tau_x : x \in \mathbb{Z}\}$ the group of translations in $\{0,1\}^{\mathbb{Z}}$ so that $(\tau_x \zeta)(z) = \zeta(x+z)$ for all $x$, $z$ in $\mathbb{Z}$ and configuration $\zeta$ in $\{0,1\}^{\mathbb{Z}}$. Translation is extended to functions and measures in a natural way.

For a positive integer $\ell$ and $-N + 1 + \ell \leq x \leq N - 1 - \ell$, denote the empirical mean density on a box of size $2\ell + 1$ centered at $x$ by $\eta^\ell(x)$,

$$\eta^\ell(x) = \frac{1}{|\Lambda_\ell(x)|} \sum_{y \in \Lambda_\ell(x)} \eta(y),$$

where $\Lambda_\ell(x) = \Lambda_{N,\ell}(x) = \{y \in \Lambda_N : |y - x| \leq \ell\}$. Let $H \in C([0,T] \times [-1,1])$ and $\Psi$ a cylinder function. For $\varepsilon > 0$ and $N$ large enough, define $V_{N,\varepsilon}^{H,\Psi} : [0,T] \times \Sigma_N \to \mathbb{R}$ by

$$V_{N,\varepsilon}^{H,\Psi}(t,\eta) = \frac{1}{N} \sum_{x=-N+1+\lfloor N\varepsilon \rfloor}^{N-1-\lfloor N\varepsilon \rfloor} H(t,x/N)\{\tau_x \Psi(\eta) - \widetilde{\Psi}(\eta^{\lfloor N\varepsilon \rfloor}(x))\},$$

where $\lfloor \cdot \rfloor$ denotes the lower integer part. Note that for $x$ as above and $N$ sufficiently large $\tau_x \Psi$ is indeed a function on $\Sigma_N$. For a function $G \in C([0,T])$ let also $W_G^\pm : [0,T] \times \Sigma_N \to \mathbb{R}$ be defined by

$$W_G^\pm(t,\eta) = G(t)[\eta(\pm N) - \rho_\pm].$$

THEOREM 3.2. *Fix $H$ in $C([0,T] \times [-1,1])$, $G$ in $C([0,T])$, a cylinder function $\Psi$, a sequence $\{\eta^N \in \Sigma_N : N \geq 1\}$ of configurations, and $\delta > 0$. Then*

$$\lim_{\varepsilon \to 0} \limsup_{N \to \infty} \frac{1}{N} \log \mathbb{P}_{\eta^N}^N \left[ \left| \int_0^T V_{N,\varepsilon}^{H,\Psi}(t,\eta_t) \, dt \right| > \delta \right] = -\infty,$$

$$\lim_{N \to \infty} \frac{1}{N} \log \mathbb{P}_{\eta^N}^N \left[ \left| \int_0^T W_G^\pm(t,\eta_t) \, dt \right| > \delta \right] = -\infty.$$

3.2. *Energy estimate.* We prove in this subsection an energy estimate. It permits the exclusion of paths with infinite energy in the large deviation regime.

Recall the definition of the constant $A_0$ introduced in Lemma 3.1. For a smooth function $G : [0,T] \times (-1,1) \to \mathbb{R}$ with compact support, let $\mathcal{Q}_G : D([0,T], \mathcal{M}) \to [0,\infty]$ be given by

$$\mathcal{Q}_G(\pi) = 2 \int_0^T dt \int_{-1}^1 du \, \pi(t,u)(\nabla G)(t,u)$$
$$- \frac{4}{A_0} \int_0^T dt \int_{-1}^1 du \, G(t,u)^2 \chi(\pi(t,u))$$

and note that

$$\mathcal{Q}(\pi) = \frac{2}{A_0} \sup_G \mathcal{Q}_G(\pi).$$

Given $\varepsilon > 0$ and a function $\pi$ in $\mathcal{M}$, let $\pi^\varepsilon : [-1,1] \to \mathbb{R}_+$ be given by

$$\pi^\varepsilon(u) = \frac{1}{2\varepsilon} \int_{[u-\varepsilon, u+\varepsilon] \cap [-1,1]} \pi(v)\, dv.$$

LEMMA 3.3. *Fix a smooth function $G : [0,T] \times (-1,1) \to \mathbb{R}$ with compact support and a sequence $\{\eta^N \in \Sigma_N : N \geq 1\}$ of configurations. There exists a constant $C > 0$, depending only on $\rho_\pm$, $E$, such that*

$$\limsup_{\varepsilon \to 0} \limsup_{N \to \infty} \frac{1}{N} \log \mathbb{P}^N_{\eta^N}[\mathcal{Q}_G(\pi^{N,\varepsilon}) \geq \ell] \leq -\ell + C(T+1).$$

PROOF. Assume without loss of generality that $\varepsilon$ is small enough for the support of $G$ to be contained in $[0,T] \times [\varepsilon, 1-\varepsilon]$. Since $\nu_N(\eta^N) \geq \exp\{-CN\}$ for some constant $C$ depending only on $\rho_\pm$, it is enough to prove the lemma with $\mathbb{P}^N_{\nu_N}$ in place of $\mathbb{P}^N_{\eta^N}$.

Let $\Psi_0(\eta) = (1/2)[\eta(1) - \eta(0)]^2$ and note that $\widetilde{\Psi}_0(\alpha) = \alpha(1-\alpha) = \chi(\alpha)$. Recall the definition of $V^{H,\Psi_0}_{N,\varepsilon}$ given just before Theorem 3.2, set $H(t,u) = G(t,u)^2$, and let $B_{N,\varepsilon}$ be the set

$$B_{N,\varepsilon} = \left\{ \eta \in D([0,T], \Sigma_N) : \left| \int_0^T V^{H,\Psi_0}_{N,\varepsilon}(t, \eta_t)\, dt \right| \leq 1 \right\}.$$

By Theorem 3.2, it is enough to show

$$\limsup_{\varepsilon \to 0} \limsup_{N \to \infty} \frac{1}{N} \log \mathbb{P}^N_{\nu_N}[\{\mathcal{Q}_G(\pi^{N,\varepsilon}) \geq \ell\} \cap B_{N,\varepsilon}] \leq -\ell + C(T+1).$$

Recall the definition of the functional $\mathcal{Q}_G$. On the one hand,

$$\int_0^T dt \int_{-1}^1 du\, (\pi^{N,\varepsilon})(t,u)(\nabla G)(t,u)$$

$$= \int_0^T dt \sum_{x=-N+1}^{N-2} \{\eta_t(x) - \eta_t(x+1)\} G(t, x/N) + O_G(\varepsilon),$$

where $O_G(\varepsilon)$ is absolutely bounded by a constant which vanishes as $\varepsilon \downarrow 0$. On the other hand, on the set $B_{N,\varepsilon}$ for $N$ large enough we have

$$\int_0^T dt \int_{-1}^1 du\, G(t,u)^2 \chi((\pi^{N,\varepsilon})(t,u))$$

$$\geq -2 + \int_0^T dt\, \frac{1}{N} \sum_{x=-N+1}^{N-2} G(t, x/N)^2 \tau_x \Psi_0(\eta_t).$$

Therefore, we just need to prove that

$$\limsup_{\varepsilon \to 0} \limsup_{N \to \infty} \frac{1}{N} \log \mathbb{P}^N_{\nu_N}\left[\int_0^T dt\, V_G(t, \eta_t) \geq \ell\right] \leq -\ell + CT$$

where

$$V_G(t, \eta) = 2 \sum_{x=-N+1}^{N-2} G(t, x/N)\{\eta(x) - \eta(x+1)\}$$

$$- \frac{2}{A_0 N} \sum_{x=-N+1}^{N-2} G(t, x/N)^2 [\eta(x) - \eta(x+1)]^2.$$

By Chebyshev's exponential inequality, the proof reduces to the statement

$$\limsup_{\varepsilon \to 0} \limsup_{N \to \infty} \frac{1}{N} \log \mathbb{E}_{\nu_N}\left[\exp\left\{N \int_0^T dt\, V_G(t, \eta_t)\right\}\right] \leq CT.$$

By Feynman–Kac's formula and the computations performed in [2], page 78, this expression is bounded by

$$\int_0^T dt \sup_f \left\{\int V_G(t, \eta) f(\eta)^2 \nu_N(d\eta) + N^{-1} \langle L_N f, f \rangle_{\nu_N}\right\},$$

where the supremum is carried over all functions $f$ in $L^2(\nu_N)$ such that $\langle f, f \rangle_{\nu_N} = 1$. By Lemma 3.1, we may replace $N^{-1} \langle L_N f, f \rangle_{\nu_N}$ by $-A_0 N D_{0,N}(f) + C_0$ for some constants $A_0, C_0$ depending only on $\rho_\pm, E$. It thus remains to show that

$$\limsup_{N \to \infty} \int_0^T dt \sup_f \left\{\int V_G(t, \eta) f(\eta)^2 \nu_N(d\eta) - A_0 N D_{0,N}(f)\right\} \leq CT$$

for some constant $C$ which depends only on $\rho_\pm$. To prove this statement we estimate the linear term of $V_G$ by its quadratic term and by $D_{0,N}(f)$.

Consider the linear term of $V_G(t, \eta)$. The change of variables $\eta' = \sigma^{x,x+1}\eta$ permits to rewrite it as

$$\sum_{x=-N+1}^{N-2} G(t, x/N) \int \{\eta(x) - \eta(x+1)\}\{f(\eta)^2 - f(\sigma^{x,x+1}\eta)^2\} \nu_N(d\eta)$$

(3.1)
$$+ \sum_{x=-N+1}^{N-2} G(t, x/N)$$

$$\times \int \{\eta(x) - \eta(x+1)\} f(\eta)^2 \{1 - F(x, \eta)\} \nu_N(d\eta),$$

where

$$F(x, \eta) = \left(\frac{\bar{\rho}(x/N)[1 - \bar{\rho}(x+1/N)]}{\bar{\rho}(x+1/N)[1 - \bar{\rho}(x/N)]}\right)^{\eta(x+1) - \eta(x)}.$$

After a Taylor expansion, the second term in (3.1) becomes

$$-\frac{1}{N}\sum_{x=-N+1}^{N-2} G(t,x/N)\nabla\left(\log\frac{\bar{\rho}}{1-\bar{\rho}}\right)(x/N)\int\{\eta(x)-\eta(x+1)\}^2 f(\eta)^2 \nu_N(d\eta)$$

plus a term of order $N^{-1}$. Since $2ab \le A_0 a^2 + A_0^{-1} b^2$, this expression is bounded by

$$C + \frac{1}{A_0 N}\sum_{x=-N+1}^{N-2} G(t,x/N)^2 \int\{\eta(x)-\eta(x+1)\}^2 f(\eta)^2 \nu_N(d\eta)$$

for some finite constant $C$ which depends on $\rho_\pm$, $E$ only. Note that the second term can be absorbed in the quadratic part of $V_G$.

To conclude the proof of the lemma, we estimate the first term in (3.1). Write $f(\eta)^2 - f(\sigma^{x,x+1}\eta)^2$ as $\{f(\eta) - f(\sigma^{x,x+1}\eta)\}\{f(\eta) + f(\sigma^{x,x+1}\eta)\}$ and apply the Schwarz inequality to bound this expression by

$$\frac{1}{4A_0 N}\sum_{x=-N+1}^{N-2} G(t,x/N)^2 \int\{\eta(x)-\eta(x+1)\}^2 \{f(\eta)+f(\sigma^{x,x+1}\eta)\}^2 \nu_N(d\eta)$$

$$+ A_0 N \sum_{x=-N+1}^{N-2} \int \{f(\eta) - f(\sigma^{x,x+1}\eta)\}^2 \nu_N(d\eta).$$

The second line is $A_0 N D_{0,N}(f)$. The first one, by a change of variables and the same arguments used to estimate the second term in (3.1), can be bounded above by

$$\frac{1}{A_0 N}\sum_{x=-N+1}^{N-2} G(t,x/N)^2 \int\{\eta(x)-\eta(x+1)\}^2 f(\eta)^2 \nu_N(d\eta) + C(G)N^{-1}$$

for some finite constant $C(G)$. This expression is part of the quadratic term of $V_G$, which concludes the proof of the lemma. □

COROLLARY 3.4. *Fix a sequence $\{G_j : j \ge 1\}$ of smooth functions $G_j : (0, T) \times (-1,1) \to \mathbb{R}$ with compact support and a sequence $\{\eta^N \in \Sigma_N : N \ge 1\}$ of configurations. There exists a constant $C$, depending only on $\rho_\pm$, $E$, such that for any $k \ge 1$*

$$\limsup_{\varepsilon \to 0}\limsup_{N\to\infty} \frac{1}{N}\log \mathbb{P}^N_{\eta^N}\left[\max_{1 \le j \le k} \mathcal{Q}_{G_j}(\pi^{N,\varepsilon}) \ge \ell\right] \le -\ell + C(T+1).$$

3.3. *Upper bound.* In this subsection, we prove the large deviations upper bound stated in Theorem 2.1. As mentioned at the beginning of this section, the proof is slightly more demanding than the original one [12, 17, 20] because the present rate function $I_T(\cdot|\gamma)$ is larger than the original one. To exclude paths with infinite energy, we rely on the estimate presented in the previous subsection.

Fix a measurable density profile $\gamma : [-1,1] \to [0,1]$, a function $H$ in $C_0^{1,2}([0,T] \times [-1,1])$, a sequence $\{G_j : j \geq 1\}$ of smooth functions $G_j : [0,T] \times (-1,1) \to \mathbb{R}$ with compact support, dense in $C_0^{0,1}([0,T] \times [-1,1])$, and a sequence of configurations $\{\eta^N : N \geq 1\}$ associated to $\gamma$. For $k \geq 1$, $\ell > 0$, let $B_{k,\ell}$ be the set of paths with truncated energy bounded by $\ell$,

$$B_{k,\ell} = \Big\{ \pi : \max_{1 \leq j \leq k} \mathcal{Q}_{G_j}(\pi) \leq \ell \Big\}.$$

By Corollary 3.4, there exists a constant $C > 0$ such that for any $k \geq 1$ and $\ell > 0$,

$$(3.2) \quad \limsup_{\varepsilon \to 0} \limsup_{N \to \infty} \frac{1}{N} \log \mathbb{P}_{\eta^N}^N [\pi^{N,\varepsilon} \notin B_{k,\ell}] \leq -\ell + C(T+1).$$

Recall the definition of the functional $\hat{J}_H : D([0,T], \mathcal{M}) \to \mathbb{R}$ introduced just before the statement of Theorem 2.1. For $k \geq 1$ and $\ell > 0$, let

$$J_H^{k,\ell}(\pi) = \begin{cases} \hat{J}_H(\pi), & \text{if } \pi \in B_{k,\ell}, \\ +\infty, & \text{otherwise.} \end{cases}$$

Let $H_1(t,u) = (\partial_u H)(t,u)$, $H_2(t,u) = H_1(t,u)^2$. Recall that $\Psi_0(\eta)$ stands for the cylinder function $(1/2)[\eta(1) - \eta(0)]^2$ and that $V_{N,\varepsilon}^{H,\Psi_0}$ is defined just before Theorem 3.2. Let $B_{\delta,\varepsilon}^{j,H,N}$, $j = 1, 2$, be the set

$$B_{\delta,\varepsilon}^{j,H,N} = \Big\{ \eta \in D([0,T], \Sigma_N) : \Big| \int_0^T V_{N,\varepsilon}^{H_j,\Psi_0}(t,\eta_t) \, dt \Big| \leq \delta \Big\}$$

and set $B_{\delta,\varepsilon}^{H,N} = B_{\delta,\varepsilon}^{1,H,N} \cap B_{\delta,\varepsilon}^{2,H,N}$. By the super-exponential estimate stated in Theorem 3.2, for each $\delta > 0$,

$$(3.3) \quad \lim_{\varepsilon \to 0} \limsup_{N \to \infty} \frac{1}{N} \log \mathbb{P}_{\eta^N}^N [(B_{\delta,\varepsilon}^{H,N})^\complement] = -\infty.$$

Fix a subset $\mathcal{A}$ of $D([0,T], \mathcal{M})$. By (3.2) and (3.3),

$$(3.4) \quad \limsup_{N \to \infty} \frac{1}{N} \log \mathbb{P}_{\eta^N}^N [\pi^N \in \mathcal{A}] \leq \max\{R_{k,\ell}^{\varepsilon,\delta}(\mathcal{A}), R_{k,\ell}(\varepsilon), R_H(\varepsilon)\},$$

where $\limsup_{\varepsilon \to 0} R_H(\varepsilon) = -\infty$, $\limsup_{\varepsilon \to 0} R_{k,\ell}(\varepsilon) \leq -\ell + C(T+1)$ and

$$R_{k,\ell}^{\varepsilon,\delta}(\mathcal{A}) = \limsup_{N \to \infty} \frac{1}{N} \log \mathbb{P}_{\eta^N}^N [\{\pi^N \in A\} \cap \{\pi^{N,\varepsilon} \in B_{k,\ell}\} \cap B_{\delta,\varepsilon}^{H,N}].$$

Consider the exponential martingale $M_t^H$ defined by

$$M_t^H = \exp\Big\{ N \Big[ \langle \pi_t^N, H_t \rangle - \langle \pi_0^N, H_0 \rangle \\ - \frac{1}{N} \int_0^t e^{-N\langle \pi_s^N, H_s \rangle} (\partial_s + L_N) e^{N\langle \pi_s^N, H_s \rangle} \, ds \Big] \Big\}.$$

Since the sequence $\{\eta^N : N \geq 1\}$ is associated to $\gamma$ and $H$ is in $C_0^{1,2}([0,T] \times [-1,1])$, an elementary computation shows that on the set $B_{\delta,\varepsilon}^{H,N}$

$$M_T^H = \exp N\{\hat{J}_H(\pi^{N,\varepsilon}) + O_H(\varepsilon) + O(\delta)\},$$

where $O_H(\varepsilon)$ [resp. $O(\delta)$] is an expression which vanishes as $\varepsilon \downarrow 0$ (resp. $\delta \downarrow 0$). On the set $\{\pi^{N,\varepsilon} \in B_{k,\ell}\}$, we may replace $\hat{J}_H(\pi^{N,\varepsilon})$ by $J_H^{k,\ell}(\pi^{N,\varepsilon})$.

Let $\mathcal{A}_{k,\ell}^{H,\varepsilon,\delta} = \{\pi^N \in \mathcal{A}\} \cap \{\pi^{N,\varepsilon} \in B_{k,\ell}\} \cap B_{\delta,\varepsilon}^{H,N}$ and write

$$\frac{1}{N} \log \mathbb{P}_{\eta^N}^N[\mathcal{A}_{k,\ell}^{H,\varepsilon,\delta}] = \frac{1}{N} \log \mathbb{E}_{\eta^N}[M_T^H (M_T^H)^{-1} \mathbf{1}\{\mathcal{A}_{k,\ell}^{H,\varepsilon,\delta}\}].$$

Optimizing over $\pi^N$ in $\mathcal{A}$, since $M_t^H$ is a mean one positive martingale, the previous expression is bounded above by

$$-\inf_{\pi \in \mathcal{A}} J_H^{k,\ell}(\pi^\varepsilon) + O_H(\varepsilon) + O(\delta).$$

Thus in view of (3.4),

$$\limsup_{N \to \infty} \frac{1}{N} \log \mathbb{P}_{\eta^N}[\pi^N \in \mathcal{A}]$$

$$\leq \max\left\{-\inf_{\pi \in \mathcal{A}} J_H^{k,\ell}(\pi^\varepsilon) + O_H(\varepsilon) + O(\delta), R_{k,\ell}(\varepsilon), R_H(\varepsilon)\right\}$$

for all $k$, $\ell$, $\varepsilon$, $\delta$ and $H$. Optimize the previous inequality with respect to these parameters and assume that the set $\mathcal{A}$ is compact. Since the map $\pi \mapsto J_H^{k,\ell}(\pi^\varepsilon)$ is lower semi-continuous for every $k$, $\ell$, $H$ and $\varepsilon$, we may apply the arguments presented in [30], Lemma 11.3, to exchange the supremum with the infimum. In this way we obtain that the last expression is bounded above by

$$\sup_{\pi \in \mathcal{A}} \inf_{H,k,\ell,\varepsilon,\delta} \max\{-J_H^{k,\ell}(\pi^\varepsilon) + O_H(\varepsilon) + O(\delta), R_{k,\ell}(\varepsilon), R_H(\varepsilon)\}.$$

For each $k \geq 1$, $\ell > 0$, and $H$ in $C_0^{1,2}([0,T] \times [-1,1])$,

$$\lim_{\varepsilon \to 0} \max_{1 \leq j \leq k} \mathcal{Q}_{G_j}(\pi^\varepsilon) = \max_{1 \leq j \leq k} \mathcal{Q}_{G_j}(\pi)$$

and $\lim_{\varepsilon \to 0} \hat{J}_H(\pi^\varepsilon) = \hat{J}_H(\pi)$. Hence $J_H^{k,\ell}(\pi) \leq \liminf_{\varepsilon \to 0} J_H^{k,\ell}(\pi^\varepsilon)$ and, letting first $\varepsilon \downarrow 0$ and then $\delta \downarrow 0$, we obtain that the previous expression is bounded above by

(3.5) $$\sup_{\pi \in \mathcal{A}} \inf_{H,k,\ell} \max\{-J_H^{k,\ell}(\pi), -\ell + C(T+1)\}.$$

Let

$$J_H^\ell(\pi) = \begin{cases} \hat{J}_H(\pi), & \text{if } \sup_{j \geq 1} \mathcal{Q}_{G_j}(\pi) \leq \ell, \\ +\infty, & \text{otherwise.} \end{cases}$$

Since the sequence $\{G_j : j \geq 1\}$ is dense in $C_0^{0,1}([0,T] \times [-1,1])$, and since $\mathcal{Q}(\pi) = (2/A_0) \sup_G \mathcal{Q}_G(\pi)$, we may replace in the previous formula the condition $\sup_{j\geq 1} \mathcal{Q}_{G_j}(\pi) \leq \ell$ by $\mathcal{Q}(\pi) \leq (2\ell/A_0)$. Since $J_H^\ell(\pi) = \lim_{k\to\infty} J_H^{k,\ell}(\pi)$, optimizing (3.5) over $k$ we obtain that it is bounded above by

$$\sup_{\pi \in \mathcal{A}} \inf_{H,\ell} \max\{-J_H^\ell(\pi), -\ell + C(T+1)\}.$$

Let

$$J_H(\pi) = \begin{cases} \hat{J}_H(\pi), & \text{if } \mathcal{Q}(\pi) < \infty, \\ +\infty, & \text{otherwise.} \end{cases}$$

Clearly, $J_H(\pi) \leq J_H^\ell(\pi)$. We may, therefore, replace in the previous variational formula $J_H^\ell(\pi)$ by $J_H(\pi)$ and let $\ell \uparrow \infty$ to conclude that the left-hand side of (3.4) is bounded above by

$$-\inf_{\pi \in \mathcal{A}} \sup_H J_H(\pi).$$

This concludes the proof of the upper bound for compact sets because $\sup_H J_H(\pi) = I_T(\pi|\gamma)$.

To pass from compact sets to closed sets, we need to prove the so-called exponential tightness for the sequence of probability measures on $D([0,T], \mathcal{M})$ given by $\{\mathbb{P}_{\eta^N}^N \circ (\pi^N)^{-1}\}$. The proof presented in [1] for the noninteracting zero-range process is easily adapted to our context.

3.4. *Lower bound.* The following is an elementary general result concerning the large deviations lower bound. Given two probability measures $P$ and $Q$ we denote by $\text{Ent}(Q|P)$ the relative entropy of $Q$ with respect to $P$.

LEMMA 3.5. *Let $\{P_n : n \geq 1\}$ be a sequence of probability measures on a Polish space $\mathcal{X}$ and let $\mathcal{X}^0 \subset \mathcal{X}$. Assume that for each $x \in \mathcal{X}^0$ there exists a sequence of probability measures $\{Q_n^x : n \geq 1\}$ which converges weakly to $\delta_x$ and such that*

(3.6) $$\limsup_n \frac{1}{n} \text{Ent}(Q_n^x|P_n) \leq I^0(x)$$

*for some functional $I^0 : \mathcal{X}^0 \to [0,\infty]$. Then, for every open set $\mathcal{O} \subset \mathcal{X}$,*

$$\liminf_{n\to\infty} \frac{1}{n} \log P_n(\mathcal{O}) \geq - \inf_{\pi \in \mathcal{O} \cap \mathcal{X}^0} I^0(x).$$

The previous result is applied with $\mathcal{X}^0$ given by the collection $\mathcal{D}_{T,\gamma}^\circ$ of "nice" paths introduced in Definition 3.6 below. In Lemma 3.7, we show that for each path $\pi$ in $\mathcal{D}_{T,\gamma}^\circ$ there exists a sequence of measures $\{\mathbb{Q}_N^\pi : N \geq 1\}$ which converges to $\pi$ and which satisfy (3.6) with $I^0(\cdot) = I_T(\cdot|\gamma)$. In view of Lemma 3.5, to complete the proof of the lower bound, it is enough to show that for every open set $\mathcal{O}$ of $D([0,T], \mathcal{M})$,

$$\text{(3.7)} \qquad \inf_{\pi \in \mathcal{O} \cap \mathcal{D}^\circ_{T,\gamma}} I_T(\pi|\gamma) = \inf_{\pi \in \mathcal{O}} I_T(\pi|\gamma).$$

This is the content of Theorem 5.1 where we prove that for any path $\pi$ in $D([0,T], \mathcal{M})$ such that $I_T(\pi|\gamma) < \infty$, there exists a sequence $\{\pi^n : n \geq 1\}$ in $\mathcal{D}^\circ_{T,\gamma}$ such that

$$\pi^n \longrightarrow \pi, \qquad I_T(\pi^n|\gamma) \longrightarrow I_T(\pi|\gamma).$$

This last property is the main subject of Section 5.

We start introducing the class of "nice" paths. This set has to be large enough to meet condition (3.7), but cannot be too large because we need to find for each nice path a sequence of probability measures satisfying the conditions of Lemma 3.5.

In the context of hydrodynamic limits, it is easy to compute the relative entropy of two dynamics which differ by a smooth external field. The external field adds to the hydrodynamic equation a transport term. The nice paths will be, therefore, paths which are solutions of the hydrodynamic equation with an extra smooth transport term. This is shown in (3.8), right after the definition of nice paths.

DEFINITION 3.6. Given $\gamma \in \mathcal{M}$, let $\mathcal{D}^\circ_{T,\gamma}$ be the collection of all paths $\rho$ in $D([0,T], \mathcal{M})$ such that:

- For every $0 < \delta \leq T$, there exists $\varepsilon > 0$ such that $\varepsilon \leq \rho(t,u) \leq 1 - \varepsilon$ for all $(t,u)$ in $[\delta, T] \times [-1, 1]$.
- There exists $\mathfrak{t} > 0$, denoted by $\mathfrak{t}(\rho)$, such that $\rho$ follows the hydrodynamic equation (2.1) in the time interval $[0, \mathfrak{t}]$, is continuous on $(0,T] \times [-1, 1]$ and smooth on $[\mathfrak{t}, T] \times [-1, 1]$.
- $\rho(0, \cdot) = \gamma(\cdot)$, $\rho(\cdot, \pm 1) = \rho_\pm$, $0 < t \leq T$.

Fix a trajectory $\pi$ in $\mathcal{D}^\circ_{T,\gamma}$. For each $0 \leq t \leq T$, let $H_t$ be the unique solution of the elliptic equation

$$\text{(3.8)} \qquad \begin{cases} \partial_t \pi_t = (1/2)\Delta \pi_t - \nabla\{\chi(\pi_t)[(E/2) + \nabla H_t]\}, \\ H_t(\pm 1) = 0. \end{cases}$$

For $t = \mathfrak{t}(\pi)$, $\partial_t \pi_t$ should be interpreted as the right derivative $\partial_{t+} \pi_t$. Note that $H$ vanishes on $[0, \mathfrak{t}(\pi)) \times [-1, 1]$. We prove in Lemma 5.7 that $H$ is smooth on $(\mathfrak{t}(\pi), T] \times [-1, 1]$ and that

$$I_T(\pi|\gamma) = \langle \pi_T, H_T \rangle - \langle \pi_{\mathfrak{t}(\pi)}, H_{\mathfrak{t}(\pi)} \rangle - \int_{\mathfrak{t}(\pi)}^T \langle \pi_t, \partial_t H_t \rangle \, dt + \frac{1}{2} \int_{\mathfrak{t}(\pi)}^T \langle \nabla \pi_t, \nabla H_t \rangle \, dt$$

$$- \frac{E}{2} \int_{\mathfrak{t}(\pi)}^T \langle \chi(\pi_t), \nabla H_t \rangle \, dt - \frac{1}{2} \int_{\mathfrak{t}(\pi)}^T \langle \chi(\pi_t), (\nabla H_t)^2 \rangle \, dt.$$

For a configuration $\eta$ in $\Sigma_N$ and a function $H:[0,T]\times[-1,1]\to\mathbb{R}$, smooth in space and smooth by parts in time, as it is the case of the function introduced in (3.8), denote by $\mathbb{P}_\eta^{N,H}$ the probability measure on $D([0,T],\Sigma_N)$ corresponding to the boundary driven weakly asymmetric exclusion process with the (weak) external field $(E/2+\nabla H)/N$ starting from $\eta$. In view of the super-exponential estimate stated in Theorem 3.2, of Definition 3.6, and of the previous explicit formula for $I_T(\pi|\gamma)$, the proof of the following lemma is similar to the one for the symmetric simple exclusion process on the torus. We thus refer to [17, 20] for its proof.

LEMMA 3.7. *Fix $\gamma \in \mathcal{M}$, a sequence $\{\eta^N \in \Sigma_N, N \geq 1\}$ associated to $\gamma$ and $\pi \in \mathcal{D}_{T,\gamma}^\circ$. Let $H$ be the solution of (3.8). The sequence of probability measures $\{\mathbb{P}_{\eta^N}^{N,H} \circ (\pi^N)^{-1}\}$ converges weakly to $\delta_\pi$ and*

$$\lim_{N\to\infty} \frac{1}{N}\operatorname{Ent}(\mathbb{P}_{\eta^N}^{N,H}|\mathbb{P}_{\eta^N}^N) = I_T(\pi|\gamma).$$

**4. The rate function.** In this and in next section we prove some properties of the rate functional $I_T$ of the large deviations principle. Since the arguments apply to a large class of interacting particle systems and might be of wider interest, we assume that the underlying stochastic dynamics has a hydrodynamical description characterized by a diffusivity $D$ and a mobility $\chi$. The method requires the mobility $\chi:[0,1]\to\mathbb{R}$ to be a smooth function equivalent to $\chi_0(\pi)=\pi(1-\pi)$ in the sense that

$$(4.1)\qquad C_0^{-1}\chi_0(\cdot) \leq \chi(\cdot) \leq C_0\chi_0(\cdot)$$

for some finite constant $C_0$, and the diffusivity $D:[0,1]\to\mathbb{R}$ to be a strictly positive smooth function. Such bounds have been proven [29, 31] for stochastic lattice gases with compact single-spin state space in the high-temperature region. We mention, however, that there are other models, such as the so-called KMP process [19], for which they do not hold.

With the previous notation the hydrodynamic equation becomes

$$(4.2)\qquad \begin{cases} \partial_t \rho = \nabla(D(\rho)\nabla\rho) - (E/2)\nabla\chi(\rho), \\ \rho(t,\pm 1) = \rho_\pm, \\ \rho(0,\cdot) = \gamma(\cdot). \end{cases}$$

For the exclusion process introduced in Section 2, $D(\rho)=1/2$ and $\chi(\rho)=\chi_0(\rho)$.

Fix once for all $T>0$. Let $\Omega=(-1,1)$ and $\Omega_T=(0,T)\times\Omega$. For a subset $E$ of $\mathbb{R}^d$, denote by $\overline{E}$ its closure. For $0\leq m, n\leq\infty$, and $E\subset\mathbb{R}$ (resp. $E\subset[0,T]\times\mathbb{R}$), denote by $C^n(E)$ [resp. $C^{m,n}(E)$] the space of functions $H:E\to\mathbb{R}$ with $n$ continuous derivatives (resp. $m$ continuous derivatives in time and $n$ continuous derivatives in space). Adding the subindex 0 (resp.

$K$) to $C^n$, $C^{m,n}$ means that the functions vanish at the boundary (resp. have compact support in the open set $E$). To keep notation simple, denote $C^{\infty,\infty}$ by $C^\infty$.

4.1. *The energy $\mathcal{Q}$.* For a *bounded* positive function $f:\Omega \to \mathbb{R}_+$ (resp. $f:\Omega_T \to \mathbb{R}_+$), denote by $L^2(f)$ [resp. $\mathbb{L}^2(f)$] the Hilbert space of (equivalence classes of) measurable functions $\{H:\Omega \to \mathbb{R}: \int_\Omega H(u)^2 f(u)\, du < \infty\}$ [resp. $\{H:\Omega_T \to \mathbb{R}: \int_{\Omega_T} H(t,u)^2 f(t,u)\, dt\, du < \infty\}$] endowed with the scalar product $\langle \cdot, \cdot \rangle_f$ (resp. $\langle\!\langle \cdot, \cdot \rangle\!\rangle_f$) induced by

$$\langle H, G \rangle_f = \int_\Omega du\, H(u) G(u) f(u), \qquad \langle\!\langle H, G \rangle\!\rangle_f = \int_0^T dt \langle H_t, G_t \rangle_{f(t,\cdot)}.$$

The norm associated to the above scalar products is denoted by $\|\cdot\|_f$. When $f = 1$ we omit the index $f$ and denote the spaces $L^2(1)$, $\mathbb{L}^2(1)$ by $L^2(\Omega)$, $L^2(\Omega_T)$, respectively.

Since $f$ is bounded, $C_K^\infty(\Omega_T)$ is dense in $\mathbb{L}^2(f)$. Moreover, the space of bounded linear functionals on $\mathbb{L}^2(f)$ can be identified with $\mathbb{L}^2(1/f)$: any bounded linear functional $\ell$ on $\mathbb{L}^2(f)$ can be represented as

(4.3) $$\ell(G) = \int_0^T dt \langle H_t, G_t \rangle$$

for some $H$ in $\mathbb{L}^2(1/f)$. Indeed, by Riesz's representation theorem, for each bounded linear functional $\ell$ on $\mathbb{L}^2(f)$, there exists a unique element $\hat{H} = \hat{H}_\ell \in \mathbb{L}^2(f)$ such that $\ell(G) = \langle\!\langle \hat{H}, G \rangle\!\rangle_f$. Let $H = f\hat{H}$. Clearly, $H$ belongs to $L^2(1/f)$ and we obtain the representation claimed above.

Fix a function $F$ in $L^1(\Omega_T)$. We claim that

$$V(F) = \sup_{G \in C_K^\infty(\Omega_T)} \left\{ 2\int_0^T dt \langle F_t, \nabla G_t \rangle - \|G\|_f^2 \right\} < \infty$$

if and only if the generalized space derivative of $F$, denoted by $\nabla F$, exists and belongs to $\mathbb{L}^2(1/f)$. In this case,

(4.4) $$V(F) = \int_0^T dt \int_{-1}^1 du\, \frac{(\nabla F_t)(u)^2}{f(t,u)}.$$

Indeed, on the one hand, assume that $V(F)$ is finite. In this case the linear functional $\ell: C_K^\infty(\Omega_T) \to \mathbb{R}$ defined by $\ell(G) = \int_0^T dt \langle F_t, \nabla G_t \rangle$ is bounded for the norm $\|\cdot\|_f$. Since $C_K^\infty(\Omega_T)$ is dense in $\mathbb{L}^2(f)$, $\ell$ can be extended to a bounded linear functional in $\mathbb{L}^2(f)$. By (4.3), for all $G$ in $C_K^\infty(\Omega_T)$,

$$\int_0^T dt \langle F_t, \nabla G_t \rangle = \ell(G) = \int_0^T dt \langle H_t, G_t \rangle$$

for some $H$ in $\mathbb{L}^2(1/f)$. Since $H$ belongs to $L^2(\Omega_T)$, this identity states that $-H$ is equal to the generalized space derivative of $F$, denoted by $\nabla F$, which belongs to $\mathbb{L}^2(1/f)$ as claimed.

Conversely, if the generalized derivative of $F$ exists and belongs to $\mathbb{L}^2(1/f)$, an integration by parts in the definition of $V(F)$ and the Schwarz inequality show that $V(f)$ is finite.

Equation (4.4) remains to be proven. After an integration by parts, the Schwarz inequality shows that the left-hand side is bounded above by the right-hand side. On the other hand, $\nabla F/f$ belongs to $\mathbb{L}^2(f)$. Since $C_K^\infty(\Omega_T)$ is dense in $\mathbb{L}^2(f)$, there exists a sequence $G_n$ in $C_K^\infty(\Omega_T)$ converging to $-\nabla F/f$ in $\mathbb{L}^2(f)$. This proves the reverse inequality.

Let $\mathcal{Q} : D([0,T], \mathcal{M}) \to [0, \infty]$ be given by

$$(4.5) \qquad \mathcal{Q}(\pi) = \frac{1}{2} \sup_{H \in C_K^\infty(\Omega_T)} \{ 2 \langle\!\langle \pi, \nabla H \rangle\!\rangle - \langle\!\langle H, H \rangle\!\rangle_{\chi_0(\pi)} \}.$$

If $\mathcal{Q}(\pi)$ is finite, by (4.4), $\pi$ has a generalized space derivative and

$$\mathcal{Q}(\pi) = \frac{1}{2} \int_0^T dt \int_{-1}^1 du \, \frac{(\nabla \pi_t)^2}{\chi_0(\pi_t)}.$$

Notice that for a path $\pi$ with $\mathcal{Q}(\pi) < \infty$, $\pi(t, \cdot)$ is continuous for almost all $t$ in $[0, T]$. While for the weakly asymmetric exclusion process $\chi = \chi_0$, in general the function $\chi$ has no reasons to be concave. It is therefore crucial to have the lower semicontinuity stated below in Lemma 4.1, that energy $\mathcal{Q}$ has been defined with $\chi_0$ and not $\chi$.

Fix a functional $J : D([0,T], \mathcal{M}) \to [0, \infty]$. A subset $\mathcal{A}$ of $D([0,T], \mathcal{M})$ is called $J$-dense if for each $\pi$ such that $J(\pi) < \infty$, there exists a sequence $\{\pi_n \in \mathcal{A} : n \geq 1\}$ converging in the topology of $D([0,T], \mathcal{M})$ to $\pi$ and such that $\lim_{n \to \infty} J(\pi_n) = J(\pi)$.

LEMMA 4.1. *The functional $\mathcal{Q}$ is convex and lower semicontinuous. The set of smooth functions bounded away from $0$ and $1$ is $\mathcal{Q}$-dense.*

PROOF. By concavity of the function $\chi_0$, for each fixed $H$ in $C_K^\infty(\Omega_T)$ the expression appearing inside braces in (4.5) is convex and therefore lower semicontinuous. These properties are inherited by $\mathcal{Q}$.

By the lower semicontinuity, to conclude the proof, it is enough to show that for each $\pi$ such that $\mathcal{Q}(\pi) < \infty$, there exists a sequence $\pi_n$ of smooth functions bounded away from 0 and 1 converging to $\pi$ and such that $\limsup_n \mathcal{Q}(\pi_n) \leq \mathcal{Q}(\pi)$.

To show that functions bounded away from 0 and 1 are $\mathcal{Q}$-dense, fix a profile $\pi$ such that $\mathcal{Q}(\pi) < \infty$. For $n \geq 1$, consider the sequence $\pi_n = n^{-1} \rho^* + (1 - n^{-1}) \pi$, where $\rho^*(t, u) = (1/2)(1-u)\rho_- + (1/2)(u+1)\rho_+$. Clearly $\pi_n$ converges to $\pi$. By convexity $\mathcal{Q}(\pi_n) \leq n^{-1}Q(\rho^*) + (1 - n^{-1})\mathcal{Q}(\pi)$. Since $\mathcal{Q}(\rho^*)$ is finite, $\limsup_n \mathcal{Q}(\pi_n) \leq \mathcal{Q}(\pi)$.

We now show that the set of smooth functions in space, bounded away from 0 and 1 are $\mathcal{Q}$-dense. Fix a function $\pi$ bounded away from 0 and 1 such that $\mathcal{Q}(\pi) < \infty$. For $\varepsilon > 0$, denote by $\lambda_\varepsilon : [-1,1] \to [-1-\varepsilon, 1+\varepsilon]$ the affine function $\lambda_\varepsilon(x) = (1+\varepsilon)x$ and by $\theta_w$ the translation by $w$ so that $(\theta_w \pi)(t,u) = \pi(t, u+w)$. Let $\pi_\varepsilon : [0,T] \times [-(1+\varepsilon), (1+\varepsilon)] \to [0,1]$ be given by $\pi_\varepsilon(t,u) = \pi(t, \lambda_\varepsilon^{-1}(u))$. We claim that

$$\mathcal{Q}(\theta_w \pi_\varepsilon) \leq \frac{1}{1+\varepsilon} \mathcal{Q}(\pi) \qquad \forall |w| \leq \varepsilon,$$

where $\theta_w \pi_\varepsilon$ must be understood as restricted to $[0,T] \times [-1,1]$. Indeed, on the one hand, for each $H$ in $C_K^\infty(\Omega_T)$,

$$\int_0^T dt \langle \theta_w \pi_\varepsilon, \nabla H \rangle = \int_0^T dt \langle \pi, \nabla J_{\varepsilon,w} \rangle,$$

where $J_{\varepsilon,w}(t,u) = H(t, \lambda_\varepsilon(u) - w)$ for $u$ in $[\lambda_\varepsilon^{-1}(-1+w), \lambda_\varepsilon^{-1}(1+w)]$ and $J_{\varepsilon,w}(t,u) = 0$ if $-1 \leq u \leq 1$ does not belong to the previous interval. On the other hand,

$$\int_0^T dt \langle \chi_0(\theta_w \pi_\varepsilon), H^2 \rangle = (1+\varepsilon) \int_0^T dt \langle \chi_0(\pi), J_{\varepsilon,w}^2 \rangle.$$

Therefore, since $J_{\varepsilon,w}$ belongs to $C_K^\infty(\Omega_T)$,

$$\int_0^T dt \langle \theta_w \pi_\varepsilon, \nabla H \rangle - \frac{1}{2} \int_0^T dt \langle \chi_0(\theta_w \pi_\varepsilon), H^2 \rangle \leq \frac{1}{1+\varepsilon} \mathcal{Q}(\pi).$$

It remains to optimize over $H$ to conclude.

Let $\{\alpha_\varepsilon : \varepsilon > 0\}$ be a smooth approximation of the identity with support contained in $(-\varepsilon, \varepsilon)$. The function $\int \alpha_\varepsilon(dw) \theta_w \pi_\varepsilon$ is smooth in space and converges to $\pi$ as $\varepsilon \downarrow 0$. By the previous estimate and the convexity of $\mathcal{Q}$,

$$\mathcal{Q}\left(\int \alpha_\varepsilon(dw) \theta_w \pi_\varepsilon\right) \leq \int \alpha_\varepsilon(dw) \mathcal{Q}(\theta_w \pi_\varepsilon) \leq \frac{1}{1+\varepsilon} \mathcal{Q}(\pi).$$

This proves that we may approximate $\pi$ by functions $\pi_n$ smooth in space and bounded away from 0 and 1 in such a way that $\mathcal{Q}(\pi_n)$ converges to $\mathcal{Q}(\pi)$.

We may repeat the same argument presented above to show that we can further require the functions to be smooth in time. $\square$

4.2. *The rate functional* $I_T(\cdot|\gamma)$. Fix once for all the initial profile $\gamma \in \mathcal{M}$. Recall that $D : [0,1] \to \mathbb{R}$ is a strictly positive continuous function and that $\chi$ is a continuous function equivalent to $\chi_0$. We next introduce the relevant rate function for interacting particle systems whose hydrodynamic behavior is described by (4.2). Let $d:[0,1] \to \mathbb{R}$ be an anti-derivative of $D$: $d' = D$, uniquely defined up to an additive constant. For each $H$ in $C_0^{1,2}(\overline{\Omega_T})$, the space of $C^{1,2}(\overline{\Omega_T})$ functions vanishing at the boundary of $\Omega$, let $\hat{J}_H = \hat{J}_{T,H,\gamma} : D([0,T], \mathcal{M}) \longrightarrow \mathbb{R}$ be the functional given by

$$
\begin{aligned}
\hat{J}_H(\pi) := {}& \langle \pi_T, H_T \rangle - \langle \gamma, H_0 \rangle - \int_0^T dt \langle \pi_t, \partial_t H_t \rangle \\
& - \int_0^T dt \langle d(\pi_t), \Delta H_t \rangle \\
& + d(\rho_+) \int_0^T dt\, \nabla H_t(1) - d(\rho_-) \int_0^T dt\, \nabla H_t(-1) \\
& - \frac{E}{2} \int_0^T dt \langle \chi(\pi_t), \nabla H_t \rangle - \frac{1}{2} \int_0^T dt \langle \chi(\pi_t), (\nabla H_t)^2 \rangle.
\end{aligned}
\tag{4.6}
$$

Of course, for the weakly asymmetric exclusion process, the above definition coincides with the one given in (2.2). Note also that the functional $\hat{J}_H$ is not affected by the choice of the arbitrary constant in the function $d$. Recalling that the general definition of the energy $\mathcal{Q}$ has been discussed in the previous subsection, the functionals $\hat{I}_T(\cdot|\gamma), I_T(\cdot|\gamma) : D([0,T], \mathcal{M}) \to [0, +\infty]$ are defined as in (2.3) and (2.4) with $\hat{J}_H$ as in (4.6). We may now state the main result of this section.

THEOREM 4.2. *For every profile $\gamma \in \mathcal{M}$, $I_T(\cdot|\gamma) : D([0,T], \mathcal{M}) \longrightarrow [0, +\infty]$ is a lower semicontinuous functional with compact level sets.*

The proof of this theorem is split in several lemmata. We show in this subsection that trajectories with finite rate function are continuous in time and satisfy the boundary conditions. We present also an alternative form of the rate function which only involves functions with compact support in $\Omega_T$.

Denote by $D_\gamma = D(\gamma, \rho_-, \rho_+)$ the subset of $D([0,T], \mathcal{M})$ of all paths $\pi$ in $C([0,T], \mathcal{M})$ which satisfy the boundary conditions $\pi(0, \cdot) = \gamma(\cdot)$, $\pi(\cdot, \pm 1) = \rho_\pm$, in the sense that the trace of $\pi$ at the boundary is $\rho_\pm$: for every $0 \le t_0 < t_1 \le T$,

$$
\lim_{\delta \to 0} \frac{1}{\delta} \int_{t_0}^{t_1} dt \int_{-1}^{-1+\delta} \pi(t,u)\, du = \rho_-(t_1 - t_0),
$$

and a similar identity at the other boundary. The proof of the next statement is similar to the one of Lemma 3.5 in [3].

LEMMA 4.3. *Fix $\pi$ in $D([0,T], \mathcal{M})$ such that $\hat{I}_T(\pi|\gamma) < \infty$. Then $\pi$ belongs to $D_\gamma$.*

In fact, for each $A > 0$, the trajectories in the set $\{\pi \in D([0,T], \mathcal{M}) : \hat{I}_T(\pi|\gamma) \leq A\}$ are uniformly continuous in time.

LEMMA 4.4. *Fix $A > 0$ and a function $J$ in $C_0^2(\overline{\Omega})$. For each $\varepsilon > 0$, there exists $\delta > 0$ such that*

$$\sup_{\pi : \hat{I}_T(\pi|\gamma) \leq A} \sup_{|t-s| \leq \delta} |\langle \pi_t, J \rangle - \langle \pi_s, J \rangle| \leq \varepsilon.$$

PROOF. Fix $A > 0$, a path $\pi$ such that $\hat{I}_T(\pi|\gamma) \leq A$ and a function $H$ in $C_0^{1,2}(\overline{\Omega_T})$. Denote by $\ell_H$ the linear part in $H$ of the functional $\hat{J}_H$. It follows from the bound $\hat{I}_T(\pi|\gamma) \leq A$ that

$$\ell_H(\pi)^2 \leq 2A \int_0^T dt \, \langle \chi(\pi_t), (\nabla H_t)^2 \rangle.$$

Fix a function $J$ in $C_0^2(\overline{\Omega})$ and $0 \leq s < t \leq T$. Approximate the indicator of the interval $[s,t]$ by smooth functions $F_\delta(r)$ and let $H^\delta(r,u) = F_\delta(r) J(u)$. With this definition,

$$\langle \pi_t, J \rangle - \langle \pi_s, J \rangle = \lim_{\delta \to 0} \left\{ \langle \pi_T, H^\delta_T \rangle - \langle \pi_0, H^\delta_0 \rangle - \int_0^T \langle \pi_r, \partial_r H^\delta_r \rangle \, dr \right\}.$$

Rewrite the expression inside braces as the sum of $\ell_{H^\delta}(\pi)$ with linear terms involving only space derivatives of $H^\delta$. Since $d, \chi$ are bounded functions, we obtain that

$$|\langle \pi_t, J \rangle - \langle \pi_s, J \rangle|$$
$$\leq C_0(t-s)\{\|J''\|_{L^1(\overline{\Omega})} + \|J'\|_{L^1(\overline{\Omega})}\} + C_0 A(t-s)^{1/2} \|J'\|_{L^2(\overline{\Omega})}$$

for some constant $C_0$ depending only on $\rho_\pm$, $E$. In this formula $\|\cdot\|_{L^p(\overline{\Omega})}$, $p \geq 1$, stand for the usual $L^p$ norm. This concludes the proof of the lemma. □

Fix $\pi$ in $D([0,T], \mathcal{M})$ such that $I_T(\pi|\gamma) < \infty$. We claim that for all $H$ in $C^{0,1}(\overline{\Omega_T})$,

$$\langle\!\langle d(\pi), \nabla H \rangle\!\rangle - d(\rho_+) \int_0^T dt \, H_t(1) + d(\rho_-) \int_0^T dt \, H_t(-1)$$
(4.7)
$$= -\langle\!\langle \nabla d(\pi), H \rangle\!\rangle,$$

where $\nabla d(\pi_t)$ stands for the generalized derivative of $d(\pi_t)$. Indeed, $\pi$ has a generalized derivative in $\mathbb{L}^2(1/\chi_0(\pi))$ because $Q(\pi) < \infty$. Thus, $d(\pi)$ has a generalized derivative which also belongs to $\mathbb{L}^2(1/\chi_0(\pi))$. Fix $H$ in $C^{0,1}(\overline{\Omega_T})$. For $\delta > 0$, let $\beta_\delta : [-1,1] \to \mathbb{R}_+$ be a smooth function with compact support in $(-1, 1)$ and equal to 1 in the interval $[-1+\delta, 1-\delta]$. Since $H_\delta(t, u) = \beta_\delta(u) H(t, u)$ belongs to $C_K^{0,1}(\Omega_T)$,

$$\langle\!\langle d(\pi), \nabla H_\delta \rangle\!\rangle = -\langle\!\langle \nabla d(\pi), H_\delta \rangle\!\rangle.$$

It remains to let $\delta \downarrow 0$ and to recall that the value of $\pi(t, \cdot)$ is fixed at the boundary to deduce (4.7).

LEMMA 4.5. *Let* $\tilde{I}_T(\cdot|\gamma) : \{\pi \in D([0, T], \mathcal{M}) : \mathcal{Q}(\pi) < \infty\} \longrightarrow [0, +\infty]$ *be the functional defined by*

$$\tilde{I}_T(\pi|\gamma) = \sup_{H \in C_K^\infty(\Omega_T)} \Big\{ -\langle\!\langle \pi, \partial_t H \rangle\!\rangle + \Big\langle\!\Big\langle D(\pi)\nabla\pi - \frac{E}{2}\chi(\pi), \nabla H \Big\rangle\!\Big\rangle \quad (4.8)$$

$$- \frac{1}{2} \langle\!\langle \chi(\pi), (\nabla H)^2 \rangle\!\rangle \Big\}.$$

*Fix* $\pi$ *in* $D_\gamma$ *such that* $\mathcal{Q}(\pi) < \infty$. *Then* $\tilde{I}_T(\pi|\gamma) = \hat{I}_T(\pi|\gamma)$.

PROOF. Fix a trajectory $\pi$ in $D_\gamma$ such that $\mathcal{Q}(\pi) < \infty$. Clearly, $\tilde{I}_T(\pi|\gamma) \leq \hat{I}_T(\pi|\gamma)$. To prove the reverse inequality, assume that $\hat{I}_T(\pi|\gamma) < \infty$ and fix $\varepsilon > 0$. By definition, there exists $H$ in $C_0^{1,2}(\Omega_T)$ such that $\hat{I}_T(\pi|\gamma) \leq \hat{J}_H(\pi) + \varepsilon$.

For $\delta > 0$, let $\beta_\delta : [0, T] \to \mathbb{R}_+$ be a smooth function with compact support in $(0, T)$ and equal to 1 in the interval $[\delta, T - \delta]$. Let

$$H_\delta(t, u) = \begin{cases} \beta_\delta(t) H(t, u/(1-\delta)), & \text{if } |u| \leq 1 - \delta, \\ 0, & \text{otherwise.} \end{cases}$$

For each $\delta > 0$, $H_\delta$ is piecewise continuously differentiable and has compact support in $\Omega_T$. Moreover, since $\pi$ belongs to $D_\gamma$, $\lim_{\delta \to 0} \hat{J}_{H_\delta}(\pi) = \hat{J}_H(\pi)$. Thus $\hat{I}_T(\pi|\gamma) \leq \hat{J}_{H_\delta}(\pi) + 2\varepsilon$ for $\delta$ small enough. It remains to approximate $H_\delta$ by a smooth function to get that $\hat{I}_T(\pi|\gamma) \leq \tilde{I}_T(\pi|\gamma)$.

On the other hand, if $\hat{I}_T(\pi|\gamma) = \infty$, one can adapt the previous arguments to show that $\tilde{I}_T(\pi|\gamma) = \infty$ as well. □

Lemmata 4.3 and 4.5 furnish an alternative definition of the rate functional $I_T(\cdot|\gamma)$:

$$(4.9) \qquad I_T(\pi|\gamma) = \begin{cases} \tilde{I}_T(\pi|\gamma), & \text{if } \pi \in D_\gamma, \mathcal{Q}(\pi) < \infty, \\ \infty, & \text{otherwise.} \end{cases}$$

We conclude this subsection with an observation on paths $\pi$ with finite energy.

LEMMA 4.6. *Fix a trajectory $\pi$ in $D_\gamma$ with finite energy: $\mathcal{Q}(\pi) < \infty$. Assume that*

$$\sup_{H \in C_K^\infty(\Omega_T)} \left\{ -\langle\!\langle \pi, \partial_t H \rangle\!\rangle - \frac{1}{4} \langle\!\langle \chi(\pi), (\nabla H)^2 \rangle\!\rangle \right\} < \infty.$$

*Then, $I_T(\pi|\gamma)$ is finite.*

PROOF. Fix a trajectory $\pi$ in $D_\gamma$ with finite energy and assume that the variational problem appearing in the statement of the lemma is finite. In view of (4.8), (4.9),

$$I_T(\pi|\gamma) \leq \sup_{H \in C_K^\infty(\Omega_T)} \left\{ -\langle\!\langle \pi, \partial_t H \rangle\!\rangle - \frac{1}{4} \langle\!\langle \chi(\pi), (\nabla H)^2 \rangle\!\rangle \right\}$$
$$+ \sup_{H \in C_K^\infty(\Omega_T)} \left\{ \left\langle\!\left\langle D(\pi) \nabla \pi - \frac{E}{2} \chi(\pi), \nabla H \right\rangle\!\right\rangle - \frac{1}{4} \langle\!\langle \chi(\pi), (\nabla H)^2 \rangle\!\rangle \right\}.$$

By assumption, the first term on the right-hand side is finite. By the Schwarz inequality, the second term is bounded above by $C_0 \mathcal{Q}(\pi)$ for some finite constant $C_0$ depending only on $E$, $D(\cdot)$ and $\chi(\cdot)$. This proves the lemma. □

4.3. *Weighted Sobolev spaces.* We introduce in this subsection weighted Sobolev spaces. We start with the classical Sobolev spaces. Let $H^1(\Omega)$ be the Sobolev space of functions $G$ in $L^2(\Omega)$ with generalized derivatives $\nabla G$ in $L^2(\Omega)$. $H^1(\Omega)$ endowed with the scalar product $\langle \cdot, \cdot \rangle_{1,2}$, defined by

$$\langle G, J \rangle_{1,2} = \langle G, J \rangle + \langle \nabla G, \nabla J \rangle,$$

is a Hilbert space. The corresponding norm is denoted by $\|\cdot\|_{1,2}$. Note that all functions in $H^1(\Omega)$ are continuous. In particular, the boundary values are well defined.

Denote by $H_0^1(\Omega)$ the closure of $C_K^\infty(\Omega)$ in $H^1(\Omega)$. Since $\Omega$ is bounded, by Poincaré's inequality, there exists a finite constant $C_1$ such that for all $G \in H_0^1(\Omega)$,

$$\|G\|_{1,2} \leq C_1 \|\nabla G\|.$$

This implies that in $H_0^1(\Omega)$,

$$\|G\|_1^2 = \langle \nabla G, \nabla G \rangle$$

is a norm equivalent to the norm $\|\cdot\|_{1,2}$. Moreover, $H_0^1(\Omega)$ is a Hilbert space with inner product given by

$$\langle G, J \rangle_1 = \langle \nabla G, \nabla J \rangle.$$

By [32], Appendix (48b), page 1030, a function $H$ in $H^1(\Omega)$ which vanishes at the boundary belongs to $H_0^1(\Omega)$.

Denote by $H^{-1}(\Omega)$ the dual of $H_0^1(\Omega)$, a Hilbert space equipped with the norm

$$\|v\|_{-1}^2 = \sup_{G \in C_K^\infty(\Omega)} \left\{ 2\langle v, G\rangle - \int_\Omega \|\nabla G(u)\|^2\, du \right\}.$$

In this formula, $\langle v, G\rangle \equiv \langle v, G\rangle_{H^{-1},H_0^1}$ stands for the value of the linear form $v$ at $G$.

Finally, for a Banach space $(\mathbb{B}, \|\cdot\|_\mathbb{B})$ and $T > 0$, we denote by $L^2(0,T;\mathbb{B})$ the Banach space of measurable functions $U:(0,T) \to \mathbb{B}$ for which

$$\|U\|_{L^2(0,T;\mathbb{B})}^2 = \int_0^T \|U(t,\cdot)\|_\mathbb{B}^2\, dt < \infty.$$

To prove the lower semicontinuity of the rate function, we need the following result which provides certain compactness.

LEMMA 4.7. *Let $\{\rho^n : n \geq 1\}$ be a sequence of functions in $L^2(\Omega_T)$ such that*

$$\int_0^T dt \|\rho_t^n\|_{1,2}^2 + \int_0^T dt \|\partial_t \rho_t^n\|_{-1}^2 \leq C_0$$

*for some finite constant $C_0$ and all $n \geq 1$. Suppose that the sequence $\rho^n$ converges weakly in $L^2(\Omega_T)$ to some $\rho$. Then $\rho^n$ converges strongly in $L^2(\Omega_T)$ to $\rho$.*

PROOF. Recall that $H^1(\Omega) \subset L^2(\Omega) \subset H^{-1}(\Omega)$. By [32], Theorem 21.A, the embedding $H^1(\Omega) \subset L^2(\Omega)$ is compact. Hence, by [28], Lemma 4, Theorem 5, the sequence $\{\rho^n : n \geq 1\}$ is relatively compact in $L^2(0,T;L^2(\Omega))$. In particular, weak convergence of the sequence $\{\rho^n : n \geq 1\}$ implies strong convergence. $\square$

We now introduce the weighted Sobolev spaces. Fix $\pi$ in $D([0,T],\mathcal{M})$ and denote by $\mathcal{H}_0^1(\chi(\pi))$ the Hilbert space induced by the smooth functions in $C_K^\infty(\Omega_T)$ endowed with the scalar product defined by

$$\langle\!\langle G, H\rangle\!\rangle_{1,\chi(\pi)} = \int_0^T dt \langle \nabla G_t, \nabla H_t\rangle_{\chi(\pi_t)}.$$

"Induced" means that we first declare two functions $F$, $G$ in $C_K^\infty(\Omega_T)$ to be equivalent if $\langle\!\langle F - G, F - G\rangle\!\rangle_{1,\chi(\pi)} = 0$ and then we complete the quotient space with respect to scalar product. Denote by $\|\cdot\|_{1,\chi(\pi)}$ the norm associated to the scalar product $\langle\!\langle \cdot, \cdot\rangle\!\rangle_{1,\chi(\pi)}$.

Let $\mathcal{H}^{-1}(\chi(\pi))$ be the dual of $\mathcal{H}_0^1(\chi(\pi))$; it is a Hilbert space equipped with the norm $\|\cdot\|_{-1,\chi(\pi)}$ defined by

$$(4.10) \qquad \|L\|_{-1,\chi(\pi)}^2 = \sup_{G \in C_K^\infty(\Omega_T)} \{2\langle\!\langle L, G\rangle\!\rangle - \|G\|_{1,\chi(\pi)}^2\}.$$

In this formula, $\langle\!\langle L, G\rangle\!\rangle$ stands for the value of the linear form $L$ at $G$. By Riesz representation theorem, an element $L$ of $\mathcal{H}^{-1}(\chi(\pi))$ can be written as $L(H) = \langle\!\langle \nabla G, \nabla H\rangle\!\rangle_{\chi(\pi)}$ for some $G$ in $\mathcal{H}_0^1(\chi(\pi))$. The next result states that $\mathcal{H}^{-1}(\chi(\pi))$ is formally the space $\{\nabla P : P \in \mathbb{L}^2(\chi(\pi)^{-1})\}$. For an integrable function $H : \Omega \to \mathbb{R}$, let $\langle H\rangle = \int_\Omega H(u)\,du$.

LEMMA 4.8. *A linear functional $L : \mathcal{H}_0^1(\chi(\pi)) \to \mathbb{R}$ belongs to $\mathcal{H}^{-1}(\chi(\pi))$ if and only if there exists $P$ in $\mathbb{L}^2(\chi(\pi)^{-1})$ such that $L(H) = \langle\!\langle P, \nabla H\rangle\!\rangle$ for every $H$ in $C_K^\infty(\Omega_T)$. In this case,*

$$\|L\|_{-1,\chi(\pi)}^2 = \int_0^T dt\{\langle P_t, P_t\rangle_{\chi(\pi_t)^{-1}} - c_t\},$$

*where $c_t = \langle P_t \chi(\pi_t)^{-1}\rangle^2 \langle \chi(\pi_t)^{-1}\rangle^{-1}\mathbf{1}\{\langle \chi(\pi_t)^{-1}\rangle < \infty\}$.*

PROOF. Fix $L$ in $\mathcal{H}^{-1}(\chi(\pi))$. By the remark preceding the lemma, $L(H) = \langle\!\langle \nabla G, \nabla H\rangle\!\rangle_{\chi(\pi)}$ for some $G$ in $\mathcal{H}_0^1(\chi(\pi))$. Let $P = \chi(\pi)\nabla G \in \mathbb{L}^2(\chi(\pi)^{-1})$ so that $L(H) = \langle\!\langle P, \nabla H\rangle\!\rangle$. Reciprocally, fix $P$ in $\mathbb{L}^2(\chi(\pi)^{-1})$. It is easy to check that the linear functional $L$ defined by $L(H) = \langle\!\langle P, \nabla H\rangle\!\rangle$ belongs to $\mathcal{H}^{-1}(\chi(\pi))$.

To compute the norm of $L$, recall that there exists $J$ in $\mathcal{H}_0^1(\chi(\pi))$ such that $L(H) = \langle\!\langle \nabla J, \nabla H\rangle\!\rangle_{\chi(\pi)}$. Therefore,

$$\int_0^T dt\langle\{P_t - \chi(\pi_t)\nabla J_t\}\nabla H_t\rangle = 0$$

for all $H$ in $C_K^\infty(\Omega_T)$. In particular, for almost all $0 \le t \le T$,

$$P_t - \chi(\pi_t)\nabla J_t = a_t \qquad \text{a.s.,}$$

where $a_t$ is a constant. The right-hand side belongs to $L^2(\chi(\pi_t)^{-1})$ because so does the left-hand side. Thus $a_t = 0$ if $\langle \chi(\pi_t)^{-1}\rangle = \infty$ and $a_t = \langle P_t\chi(\pi_t)^{-1}\rangle\langle \chi(\pi_t)^{-1}\rangle^{-1}$, otherwise. Moreover,

$$\|L\|_{-1,\chi(\pi)}^2 = \sup_{H \in C_K^\infty(\Omega_T)} \{2\langle\!\langle P, \nabla H\rangle\!\rangle - \|H\|_{1,\chi(\pi)}^2\}$$

$$= \sup_{H \in C_K^\infty(\Omega_T)} \{2\langle\!\langle \nabla J, \nabla H\rangle\!\rangle_{\chi(\pi)} - \|H\|_{1,\chi(\pi)}^2\}$$

$$= \int_0^T dt\langle \chi(\pi_t)(\nabla J_t)^2\rangle.$$

To conclude the proof, it remains to recall that $\chi(\pi_t)\nabla J_t = P_t - a_t$. □

Fix an integrable function $G:[0,T] \times [-1,1] \to \mathbb{R}$. Then,
$$V(G) = \sup_{H \in C_K^\infty(\Omega_T)} \{2\langle\!\langle G, H\rangle\!\rangle - \|H\|_{1,\chi(\pi)}^2\} < \infty$$

if and only if the linear functional $L_G: C_K^\infty(\Omega_T) \to \mathbb{R}$ defined by $L_G(H) = \langle\!\langle G, H\rangle\!\rangle$ belongs to $\mathcal{H}^{-1}(\chi(\pi))$. Indeed, if $V(G) < \infty$, $L_G$ is bounded in $\mathcal{H}_0^1(\chi(\pi))$ and thus belongs to $\mathcal{H}^{-1}(\chi(\pi))$. On the other hand, if $L_G$ belongs to $\mathcal{H}^{-1}(\chi(\pi))$, by Lemma 4.8, for each $H$ in $C_K^\infty(\Omega_T)$, $\langle\!\langle G, H\rangle\!\rangle = L_G(H) = \langle\!\langle P, \nabla H\rangle\!\rangle$ for some $P$ in $\mathbb{L}^2(\chi(\pi)^{-1})$. Hence, $V(G)$ is finite by the Schwarz inequality. In this case, $V(G) = \|L_G\|_{-1,\chi(\pi)}^2$.

By abuse of notation, we shall say that $G$ belongs to $\mathcal{H}^{-1}(\chi(\pi))$ whenever the linear functional $L_G$ belongs to $\mathcal{H}^{-1}(\chi(\pi))$. In this case we denote by $\|G\|_{-1,\chi(\pi)}$ the norm of $L_G$: $\|G\|_{-1,\chi(\pi)} = \|L_G\|_{-1,\chi(\pi)}$. With this convention, recalling (4.8), for every path $\pi$ in $D([0,T],\mathcal{M})$ with finite energy, $\mathcal{Q}(\pi) < \infty$,

(4.11)    $\tilde{I}_T(\pi|\gamma) = \tfrac{1}{2}\|\partial_t \pi - \nabla(D(\pi)\nabla\pi) + (E/2)\nabla\chi(\pi)\|_{-1,\chi(\pi)}^2$,

where this expression might take the value $+\infty$. In this formula, $\partial_t \pi$ is the linear functional whose value at $H \in C_K^\infty(\Omega_T)$ is equal to $-\langle\!\langle \pi_t, \partial_t H\rangle\!\rangle$.

Fix a trajectory $\pi$ in $D_\gamma$ with finite energy, $\mathcal{Q}(\pi) < \infty$. Since $\chi$ and $\chi_0$ are equivalent and since $D$ is bounded, the weak derivatives $\nabla(D(\pi)\nabla\pi)$, $\nabla\chi(\pi)$ belong to $\mathcal{H}^{-1}(\chi(\pi))$. In particular, $\partial_t\pi$ belongs to $\mathcal{H}^{-1}(\chi(\pi))$ if and only if $I_T(\pi|\gamma)$ is finite. Indeed, on the one hand, if $\|\partial_t\pi\|_{-1,\chi(\pi)} < \infty$, by Lemma 4.6, $I_T(\pi|\gamma) < \infty$. On the other hand, if $I_T(\pi|\gamma) < \infty$, it follows from (4.11) that $\partial_t\pi$ belongs to $\mathcal{H}^{-1}(\chi(\pi))$ as well.

If $\partial_t\pi$ belongs to $\mathcal{H}^{-1}(\chi(\pi))$, by Lemma 4.8, $\partial_t\pi = \nabla P$ for some $P = P_\pi$ in $\mathbb{L}^2(\chi(\pi)^{-1})$,

(4.12)    $\langle\!\langle \pi, \partial_t H\rangle\!\rangle = \langle\!\langle P, \nabla H\rangle\!\rangle$

for every $H$ in $C_K^\infty(\Omega_T)$. We may of course choose $P$ so that $\langle P_t \chi(\pi_t)^{-1}\rangle \times \langle \chi(\pi_t)^{-1}\rangle^{-1} \mathbf{1}\{\langle \chi(\pi_t)^{-1}\rangle < \infty\} = 0$. Replacing $\langle\!\langle \pi, \partial_t H\rangle\!\rangle$ by $\langle\!\langle P, \nabla H\rangle\!\rangle$ in the variational formula appearing in the statement of Lemma 4.5, we obtain from Lemma 4.8 an explicit expression for the rate functional,

(4.13)   $I_T(\pi|\gamma) = \dfrac{1}{2}\int_0^T dt\{\|P_t - D(\pi_t)\nabla\pi_t + (E/2)\chi(\pi_t)\|_{\chi(\pi_t)^{-1}}^2 - R_t\}$,

where
$$R_t = \{\delta h - E\}^2 \frac{1}{\langle\chi(\pi_t)^{-1}\rangle}$$

and $\delta h = h(\rho_+) - h(\rho_-)$ where $h'(\rho) = D(\rho)/\chi(\rho)$. Here we adopted the convention that $R_t$ vanishes if $\langle\chi(\pi_t)^{-1}\rangle = \infty$. Note that $R_t$ vanishes in the reversible case because $h(\rho_+) - h(\rho_-) = E$ there.

LDP FOR THE WEAKLY ASYMMETRIC EXCLUSION PROCESS 29

4.4. *Lower-semicontinuity of $I_T(\cdot|\gamma)$*. In this subsection we conclude the proof that the rate function $I_T(\cdot|\gamma)$ is lower-semicontinuous and has compact level sets.

LEMMA 4.9. *There exists a constant $C_0$ such that*

$$\|\partial_t \rho\|^2_{-1,\chi(\rho)} \leq C_0\{1 + I_T(\rho|\gamma) + \mathcal{Q}(\rho)\}, \mathcal{Q}(\rho) \leq C_0\{1 + I_T(\rho|\gamma)\}$$

*for all $\rho$ in $D([0,T], \mathcal{M})$.*

PROOF. Fix $\rho$ in $D([0,T], \mathcal{M})$ such that $I_T(\rho|\gamma) < \infty$. By (4.11),

$$\|\partial_t \rho\|^2_{-1,\chi(\rho)} \leq 4 I_T(\rho|\gamma) + 2\|\nabla(D(\rho)\nabla\rho) - (E/2)\nabla\chi(\rho)\|^2_{-1,\chi(\rho)}.$$

Recall that $D$ is bounded. By definition of the norm $\|\cdot\|_{-1,\chi(\rho)}$ and by Schwarz inequality, the second term is bounded above by $C_0\{1 + \mathcal{Q}(\rho)\}$ for some finite constant $C_0$ which depends only on $D$, $E$ and $\chi$. This concludes the first part of the proof.

We now prove the second statement of the lemma. Since $\mathcal{Q}(\rho) < \infty$ it follows that $\rho(t,\cdot)$ is continuous and $\rho(t,\pm 1) = \rho_\pm$ for almost all $t \in [0,T]$, and $\partial_t \rho$ belongs to $L^2([0,T]; H^{-1}(\Omega))$.

Let $\varepsilon = \min\{\rho_-, 1 - \rho_+\} > 0$ and let $\mathcal{D} = \{(a,b) \in [0,1] \times [-1,1] : 0 \leq a + b \leq 1\}$. Fix $0 < \delta < \varepsilon$ small and let $F_\delta : \mathcal{D} \to \mathbb{R}$ be defined by

$$F_\delta(a,b) = (a + b + \delta) \log \frac{a + b + \delta}{a + \delta} + (1 + \delta - a - b) \log \frac{1 + \delta - a - b}{1 + \delta - a}.$$

Fix a smooth profile $\bar{\rho} \in \mathcal{M}$ satisfying the boundary conditions, $\bar{\rho}(\pm 1) = \rho_\pm$. Since $\rho$ belongs to $L^2([0,T]; H^1(\Omega))$ and $\rho(t,\pm 1) = \rho_\pm$ for almost all $t \in [0,T]$, by [32], Appendix (48b), page 1030, $\rho - \bar{\rho}$ belongs to $L^2([0,T]; H^1_0(\Omega))$. As $\partial_t(\rho - \bar{\rho}) = \partial_t \rho$ belongs to $L^2([0,T]; H^{-1}(\Omega))$, by [32], Proposition 23.23(iii), there exists a sequence $\{U^n : n \geq 1\}$ of smooth functions $U^n : [0,T] \times \Omega \to \mathbb{R}$ with compact support such that $U^n$ converges to $\rho - \bar{\rho}$ in $L^2([0,T]; H^1_0(\Omega))$ and $\partial_t U^n$ converges to $\partial_t \rho$ in $L^2([0,T]; H^{-1}(\Omega))$.

For each $n \geq 1$,

$$\int_\Omega F_\delta(\bar{\rho}, U^n_T) - \int_\Omega F_\delta(\bar{\rho}, U^n_0) = \int_0^T dt \langle \partial_t U^n_t, (\partial_2 F_\delta)(\bar{\rho}, U^n_t)\rangle,$$

where $\partial_2 F_\delta$ stands for the partial derivative of $F_\delta$ with respect to the second coordinate. Since $F_\delta$ is smooth with bounded first and second derivatives and since $U^n$ converges to $\rho - \bar{\rho}$, letting $n \uparrow \infty$, we obtain that

$$0 = \int_\Omega F_\delta(\bar{\rho}, \rho_T - \bar{\rho}) - \int_\Omega F_\delta(\bar{\rho}, \rho_0 - \bar{\rho}) - \langle\!\langle \partial_t \rho, (\partial_2 F_\delta)(\bar{\rho}, \rho - \bar{\rho})\rangle\!\rangle.$$

In this formula, the scalar product on the right-hand side has to be understood as the value at $(\partial_2 F_\delta)(\bar{\rho}, \rho - \bar{\rho})$ of the linear functional $\partial_t \rho$. Since the last term is equal to $\langle\!\langle \partial_t \rho, h_\delta(\rho) - h_\delta(\bar\rho) \rangle\!\rangle$ where $h_\delta(x) = \log\{\delta + x/1 + \delta - x\}$, the previous identity can be written as

$$(4.14) \quad \langle\!\langle \partial_t \rho, h_\delta(\rho_t) - h_\delta(\bar\rho) \rangle\!\rangle = \int_\Omega F_\delta(\bar\rho, \rho_T - \bar\rho) - \int_\Omega F_\delta(\bar\rho, \rho_0 - \bar\rho).$$

Note that the right-hand side is absolutely bounded uniformly in $\delta$, by a constant depending only on $\bar\rho$.

Fix a path $\pi$ in $\mathcal{H}^1(\chi(\pi))$. We claim that for any $\alpha > 0$ and any $H$ in $L^2([0,T]; H_0^1(\Omega))$.

$$(4.15) \quad \langle\!\langle \partial_t \pi, H \rangle\!\rangle \leq \frac{1}{\alpha} I_T(\pi|\gamma) - \langle\!\langle D(\pi)\nabla\pi, \nabla H \rangle\!\rangle + \frac{E}{2}\langle\!\langle \chi(\pi_t), \nabla H_t \rangle\!\rangle + \frac{\alpha}{2}\|H\|^2_{1,\chi(\pi)}.$$

We prove this statement at the end of the lemma.

Since $h_\delta(\rho_t) - h_\delta(\bar\rho)$ belongs to $L^2([0,T]; H_0^1(\Omega))$, (4.15) holds with $\pi$ replaced by $\rho$ and with $H_t = h_\delta(\rho_t) - h_\delta(\bar\rho)$. By the Schwarz inequality, the third term on the right-hand side is bounded by $C_0 \alpha^{-1} + (\alpha/2)\|H\|^2_{1,\chi(\pi)}$ for every $\alpha > 0$. Hereafter, $C_0$ stands for a finite constant, depending only on $E$, $\rho_\pm$, $D$ and $\chi$, whose value may change from line to line. We have seen just after (4.14) that the left-hand side is absolutely bounded by a constant depending only on $\bar\rho$. Hence, moving the term on the right-hand side to the left-hand side and the second term on the right-hand side to the left-hand side, we get that

$$\langle\!\langle D(\rho_t)\nabla\rho_t, \nabla H_t \rangle\!\rangle \leq \frac{1}{\alpha} C_0 + \frac{1}{\alpha} I_T(\rho|\gamma) + \alpha \|H\|^2_{1,\chi(\pi)}$$

for every $\alpha > 0$. Let $\chi_\delta(a) = (a + \delta)(1 + \delta - a)/(1 + 2\delta)$ so that $h'_\delta(a) = \chi_\delta(a)^{-1}$. Since $H_t = h_\delta(\rho_t) - h_\delta(\bar\rho)$ and $D$ is bounded below by a strictly positive constant, we may choose $\alpha$ small enough to get

$$\int_0^T dt \int_\Omega du \frac{(\nabla\rho_t)^2}{\chi_\delta(\rho_t)} \leq C_0\{1 + I_T(\rho|\gamma)\}$$

for some constant $C_0$. Applying Fatou's lemma we obtain that $\mathcal{Q}(\rho) \leq C_0\{1 + I_T(\rho|\gamma)\}$.

It remains to prove (4.15). Fix a path $\pi$ in $\mathcal{H}^1(\chi(\pi))$ and $H$ in $C_K^\infty(\Omega_T)$. By the explicit formula (4.11) for $I_T(\pi|\gamma)$ and the variational formula (4.10) for the norm $\|\cdot\|_{-1,\chi(\pi)}$,

$$\langle\!\langle \partial_t \pi - \nabla(D(\pi)\nabla\pi) + (E/2)\nabla\chi(\pi), H \rangle\!\rangle - \frac{\alpha}{2}\|H\|^2_{1,\chi(\pi)} \leq \frac{1}{\alpha} I_T(\pi|\gamma)$$

for every $\alpha > 0$. To conclude the proof, it remains to recall that $C_K^\infty([0,T] \times \Omega)$ is dense in $L^2([0,T]; H_0^1(\Omega))$. $\square$

PROOF OF THEOREM 4.2. To prove the lower semicontinuity we have to show that for all $\lambda \geq 0$, the set
$$E_\lambda = \{\pi \in D([0,T], \mathcal{M}) : I_T(\pi|\gamma) \leq \lambda\}$$
is closed in $D([0,T], \mathcal{M})$. Fix $\lambda \geq 0$ and consider a sequence $\{\pi_n : n \geq 1\}$ in $E_\lambda$ converging to some $\pi$ in $D([0,T], \mathcal{M})$. Thus for all $G$ in $C(\overline{\Omega_T})$,

(4.16) $$\lim_{n\to\infty} \langle\!\langle G, \pi_n \rangle\!\rangle = \langle\!\langle G, \pi \rangle\!\rangle.$$

By Lemma 4.9, there exists a positive constant $C_\lambda$ such that
$$\sup_{n\geq 1} \mathcal{Q}(\pi_n) \leq C_\lambda \quad \text{and} \quad \sup_{n\geq 1} \int_0^T dt \|\partial_t \pi_n\|_{-1}^2 \leq C_\lambda.$$

By (4.16), $\pi_n$ converges weakly to $\pi$ in $L^2(\Omega_T)$. Hence by Lemma 4.7, $\pi_n$ converges strongly to $\pi$ in $L^2(\Omega_T)$. The proof of Lemma 4.1 shows that $\mathcal{Q}$ is lower semicontinuous also for the strong $L^2(\Omega_T)$ topology so that $\mathcal{Q}(\pi) \leq C_\lambda < \infty$.

Fix $G$ in $C_K^\infty(\Omega_T)$. Since $\pi_n$ converges strongly to $\pi$ in $L^2(\Omega_T)$,
$$\lim_{n\to\infty} \{\langle\!\langle \pi_n, \partial_t G \rangle\!\rangle - \langle\!\langle D(\pi_n)\nabla\pi_n, \nabla G \rangle\!\rangle + (E/2)\langle\!\langle \chi(\pi_n), \nabla G \rangle\!\rangle - \|G\|_{1,\chi(\pi_n)}^2\}$$
$$= \langle\!\langle \pi, \partial_t G \rangle\!\rangle - \langle\!\langle D(\pi)\nabla\pi, \nabla G \rangle\!\rangle + (E/2)\langle\!\langle \chi(\pi), \nabla G \rangle\!\rangle - \|G\|_{1,\chi(\pi)}^2.$$

Since $\pi_n$ belongs to $E_\lambda$, the left-hand side is bounded by $\lambda$. Taking the supremum over $G$ in $C_K^\infty(\Omega_T)$ we obtain that $\tilde{I}_T(\pi|\gamma) \leq \lambda$.

We claim that $\pi$ belongs to $D_\gamma$. The proof of Lemma 4.4 with $\tilde{I}_T(\cdot|\gamma)$ in place of $\hat{I}_T(\cdot|\gamma)$ shows that $\pi$ is uniformly continuous in time. In particular, $\pi$ belongs to $C([0,T], \mathcal{M})$. Furthermore, since $\pi_n \in D_\gamma$ converges to $\pi$ in $L^2(\Omega_T)$ and $\pi \in C([0,T], \mathcal{M})$, $\pi(0, \cdot) = \gamma(\cdot)$.

To show that $\pi(\cdot, \pm 1) = \rho_\pm$, recall that the boundary values are well defined for any function $\rho$ in $H^1(\Omega)$. Moreover,
$$\int \nabla \rho H \, du = H(1)\rho(1) - H(-1)\rho(-1) - \int \rho \nabla H \, du$$
for all function $H$ in $C^\infty([-1,1])$. Since $\mathcal{Q}(\pi_n) \leq C_\lambda$, there exists a subsequence $n_k$ and $v$ in $L^2(\Omega_T)$ such that $\nabla \pi_{n_k}$ converges weakly in $L^2(\Omega_T)$ to $v$. Since $\pi_n$ converges in $L^2(\Omega_T)$ to $\pi$, by an integration by parts formula for time-dependent smooth functions $H$ with compact support in $[0,T] \times \Omega$, similar to the last displayed equation, $v_t = \nabla \pi_t$ for almost all $t$. Since $\pi_n$ belongs to $D_\gamma$, again by the integration by parts formula $\pi(\cdot, \pm 1) = \rho_\pm$ for almost all $t$. This proves the claim that $\pi$ belongs to $D_\gamma$.

To conclude the proof of the lower semicontinuity, note that $I_T(\pi|\gamma) \leq \lambda$ in view of (4.9) and the estimates $\mathcal{Q}(\pi) < \infty$, $\tilde{I}_T(\pi|\gamma) \leq \lambda$ obtained above.

We now turn to the proof of the compact level sets. Consider a sequence of trajectories $\{\rho_n : n \geq 1\}$ such that $I_T(\rho_n|\gamma) \leq \lambda$. Since each trajectory is positive and bounded by 1, there exists a subsequence, still denoted by $\{\rho_n : n \geq 1\}$, which converges weakly in $L^2(\Omega_T)$ to some trajectory $\rho$. Repeating the arguments presented in the first part of the proof, we may conclude that $\rho_n$ converges strongly to $\rho$ in $L^2(\Omega_T)$ and that $\mathcal{Q}(\rho) < \infty$. The first part of the proof shows also that $I_T(\cdot|\gamma)$ is lower semicontinuous for the weak $L^2(\Omega_T)$ topology so that $I_T(\rho|\gamma) \leq \lambda$. By Lemma 4.4, $\rho_n$, $\rho$ are uniformly continuous in time. In particular, strong convergence in $L^2(\Omega_T)$ implies convergence in $C([0,T], \mathcal{M})$. $\square$

**5. $I_T(\cdot|\gamma)$-density.** In this section we show that any path $\pi \in D([0,T], \mathcal{M})$ with finite rate function, $I_T(\pi|\gamma) < \infty$, can be approximated by the smooth paths introduced in Definition 3.6. As in the previous section, we work with an arbitrary smooth diffusion coefficient $D$ uniformly positive and an arbitrary mobility $\chi$ which satisfies the bounds (4.1). In particular, in the Definition 3.6 we need to replace the hydrodynamic equation (2.1) by (4.2). The main theorem of this section is stated as follows.

THEOREM 5.1. *Fix $\gamma \in \mathcal{M}$. The set $\mathcal{D}_{T,\gamma}^\circ$ is $I_T(\cdot|\gamma)$-dense.*

The proof of the $I_T(\cdot|\gamma)$-density of some set $\mathcal{A}$ relies on the next two results. Recall from (4.12) that for each path $\pi$ such that $I_T(\pi|\gamma) < \infty$, there exists $P = P_\pi$ in $\mathbb{L}^2(\chi(\pi)^{-1})$ such that $\langle\!\langle \pi, \partial_t H \rangle\!\rangle = \langle\!\langle P, \nabla H \rangle\!\rangle$ for every $H$ in $C_K^\infty(\Omega_T)$ and $\langle P_t \chi(\pi_t)^{-1} \rangle \mathbf{1}\{\langle \chi(\pi_t)^{-1} \rangle < \infty\} = 0$.

LEMMA 5.2. *Fix a trajectory $\pi$ with $I_T(\pi|\gamma) < \infty$. Consider a sequence $\{\pi_n : n \geq 1\}$ such that $I_T(\pi_n|\gamma) < \infty$ and:*

1. *$\pi_n$, $P_n = P_{\pi_n}$, $\nabla \pi_n$ converge to $\pi$, $P_\pi$, $\nabla \pi$ almost everywhere in $\Omega_T$;*
2. *$\{D(\pi_n)\nabla \pi_n\}^2/\chi(\pi_n)$, $P_{\pi_n}^2/\chi(\pi_n)$ are uniformly integrable;*
3. *$\int_0^T dt\, 1/\langle \chi(\pi_n(t))^{-1} \rangle$ converges to $\int_0^T dt\, 1/\langle \chi(\pi(t))^{-1} \rangle$.*

*Then, $I_T(\pi_n|\gamma)$ converges to $I_T(\pi|\gamma)$.*

PROOF. Recall that $D(\cdot)$, $\chi(\cdot)$ are continuous functions and recall the explicit form (4.13) for the rate function $I_T(\pi|\gamma)$. The assumptions of the lemma are tailored for $I_T(\pi_n|\gamma)$ to converge to $I_T(\pi|\gamma)$. $\square$

The following elementary lemma will be used repeatedly in the sequel to prove uniform integrability of sequences of functions.

LEMMA 5.3. *Fix a measure space $(\Omega, \mu, \mathcal{F})$ and a function $f$ in $L^1(\mu)$. There exists an increasing convex function $\Psi : \mathbb{R}_+ \to \mathbb{R}_+$ such that $\lim_{x \to \infty} \Psi(x)/x = \infty$ and*

$$\int \Psi(|f|)\, d\mu < \infty.$$

*A family $\{f_\alpha\}$ of functions satisfying*

$$\sup_\alpha \int \Phi(|f_\alpha|)\, d\mu < \infty$$

*for a function $\Phi$ such that $\lim_{x\to\infty} \Phi(x)/x = \infty$ is uniformly integrable.*

PROOF. For $x \geq 0$ let $G(x) := \int_x^\infty dy\, \mu(|f| > y)$. Then $G(0) = \int d\mu\, |f| < \infty$ and $G(x) \downarrow 0$ as $x \uparrow \infty$. It is simple to check that the function $\Psi(x) = \int_0^x dy\, G(y)^{-1/2}$ meets the requirements of the lemma. The second statement is trivial. $\square$

The proof of Theorem 5.1 strongly uses the smoothing effect of the hydrodynamic equation. Denote by $\rho$ the solution of the hydrodynamic equation (4.2) with initial condition $\gamma$ so that $\hat{I}_T(\rho|\gamma) = 0$. We claim that

$$(5.1) \qquad \int_0^T dt \int_{-1}^1 du\, \frac{\{\nabla \rho\}^2}{\chi(\rho)} < \infty, \qquad \int_0^T dt \int_{-1}^1 du\, \frac{P_\rho^2}{\chi(\rho)} < \infty.$$

Let $F:[0,1] \to \mathbb{R}_+$ such that $F''(\alpha) = D(\alpha)/\chi(\alpha)$, $\alpha \in (0,1)$. In view of our assumptions, $F$ is bounded but its derivative $F'(\alpha)$ diverges as $\alpha \to 0, 1$. Pick a sequence of smooth functions $F_n: \mathbb{R} \to \mathbb{R}_+$ such that $F_n(\alpha) \uparrow F(\alpha)$, $F_n''(\alpha) \uparrow F''(\alpha)$, $\alpha \in (0,1)$, as $n \to \infty$. The first estimate in (5.1) is proven by computing the time derivative of $\int_\Omega F_n(\rho(t,u))\, du$ and taking the limit $n \to \infty$.

Since, by (4.1), $\chi$ is equivalent to $\chi_0$, the first estimate in (5.1) shows that the energy of $\rho$ is finite: $\mathcal{Q}(\rho) < \infty$. Therefore, by definition of $I_T(\cdot|\gamma)$, $I_T(\rho|\gamma) = \hat{I}_T(\rho|\gamma)$. Since $\hat{I}_T(\rho|\gamma)$ vanishes, $I_T(\rho|\gamma) = 0$. In particular, from the explicit formula (4.13) for $I_T(\rho|\gamma)$ we have that

$$\int_0^T dt \|P_t\|_{\chi(\rho_t)^{-1}}^2 \leq 2 \int_0^T dt \|D(\rho_t)\nabla \rho_t - (E/2)\chi(\rho_t)\|_{\chi(\rho_t)^{-1}}^2 + 2 \int_0^T dt\, R_t.$$

The finiteness of this expression follows from the boundedness of $\chi$ and from the first estimate in (5.1).

We are now ready to prove the first result towards Theorem 5.1. Let $\mathcal{F}_0$ be the subset of $D_\gamma$ of all trajectories $\pi$ such that $I_T(\pi|\gamma) < \infty$ and for which there exist $\delta > 0$ such that $\pi$ is equal to the solution of the hydrodynamic equation (4.2) in the time interval $[0, \delta]$. More precisely, denote by $\rho$ the solution of (4.2). There exists $\delta > 0$ such that $\pi(t,u) = \rho(t,u)$ for $(t,u)$ in $[0, \delta] \times [-1,1]$.

LEMMA 5.4. *The set $\mathcal{F}_0$ is $I_T(\cdot|\gamma)$-dense.*

PROOF. Fix a path $\pi$ such that $I_T(\pi|\gamma) < \infty$ and let $\rho$ be the solution of the hydrodynamic equation (4.2). For $\varepsilon > 0$, define $\pi^\varepsilon$ as

$$\pi^\varepsilon(t,\cdot) = \begin{cases} \rho(t,\cdot), & \text{for } 0 \leq t \leq \varepsilon, \\ \rho(2\varepsilon - t, \cdot), & \text{for } \varepsilon \leq t \leq 2\varepsilon, \\ \pi(t - 2\varepsilon, \cdot), & \text{for } 2\varepsilon \leq t \leq T. \end{cases}$$

For each $\varepsilon > 0$, $\pi^\varepsilon$ belongs to $D_\gamma$ because so does $\pi$ and because $\rho$ is the solution of the hydrodynamic equation. Moreover, $\mathcal{Q}(\pi^\varepsilon) \leq \mathcal{Q}(\pi) + 2\mathcal{Q}(\rho) < \infty$ and $\pi^\varepsilon$ converges to $\pi$ as $\varepsilon \downarrow 0$ because $\pi$ belongs to $C([0,T], \mathcal{M})$. It remains to show that $\hat{I}_T(\pi^\varepsilon|\gamma)$ converges to $\hat{I}_T(\pi|\gamma)$.

By lower semicontinuity, $\hat{I}_T(\pi|\gamma) \leq \liminf_{\varepsilon \to 0} \hat{I}_T(\pi^\varepsilon|\gamma)$. To prove the reverse inequality, decompose the rate function $\hat{I}_T(\pi^\varepsilon|\gamma)$ as the sum of the contributions on each time interval $[0, \varepsilon]$, $[\varepsilon, 2\varepsilon]$ and $[2\varepsilon, T]$. The first contribution vanishes because $\pi^\varepsilon$ follows the hydrodynamic equation in this interval and the third contribution is bounded by $\hat{I}_T(\pi|\gamma)$ because $\pi^\varepsilon$ in this interval is just a time translation of the path $\pi$.

On the interval $[\varepsilon, 2\varepsilon]$, $\pi^\varepsilon$ is the solution of the equation

$$\partial_t \rho_t = -\nabla(D(\rho)\nabla\rho) + (E/2)\nabla\chi(\rho).$$

In particular, by Lemma 4.5, the contribution of the interval $[\varepsilon, 2\varepsilon]$ to the rate function is equal to

$$\sup_{H \in C_K^\infty(\Omega_T)} \left\{ 2\int_0^\varepsilon dt \left\langle D(\rho_t)\nabla\rho_t - \frac{E}{2}\chi(\rho_t), \nabla H_t \right\rangle - \frac{1}{2}\int_0^\varepsilon dt \langle \chi(\rho_t), (\nabla H_t)^2 \rangle \right\}.$$

By the Schwarz inequality, the previous expression is less than or equal to

$$C_0 \left\{ \varepsilon + \int_0^\varepsilon dt \left\langle \frac{(\nabla \rho_t)^2}{\chi(\rho_t)} \right\rangle \right\}$$

for some finite constant $C_0$ which depends only on $E$, $D(\cdot)$ and $\chi(\cdot)$. Since, by (5.1), $\mathcal{Q}(\rho)$ is finite, this expression vanishes as $\varepsilon \downarrow 0$ and we are done. $\square$

Denote by $\mathcal{F}_1$ the subset of $\mathcal{F}_0$ of all trajectories $\pi$ such that for all $0 < \delta \leq T$, there exists $\varepsilon > 0$ such that $\varepsilon \leq \pi(t, u) \leq 1 - \varepsilon$ for $\delta \leq t \leq T$, $-1 \leq u \leq 1$.

LEMMA 5.5. *The set $\mathcal{F}_1$ is $I_T(\cdot|\gamma)$-dense.*

PROOF. Fix $\pi$ in $\mathcal{F}_0$. In view of the previous lemma, it is enough to exhibit a sequence $\{\pi_\varepsilon : \varepsilon > 0\}$ in $\mathcal{F}_1$ which converges to $\pi$ and such that $I_T(\pi_\varepsilon|\gamma)$ converges to $I_T(\pi|\gamma)$.

For $0 < \varepsilon < 1$, let $\pi_\varepsilon = (1-\varepsilon)\pi + \varepsilon\rho$ where $\rho$ is the solution of the hydrodynamic equation (4.2). We claim that $\pi_\varepsilon$ belongs to $\mathcal{F}_1$ for each $0 < \varepsilon < 1$.

On the one hand, $\pi_\varepsilon$ belongs to $D_\gamma$ because so do $\pi$ and $\rho$. Moreover, $\mathcal{Q}(\pi_\varepsilon) \leq C_0 < \infty$ because $\mathcal{Q}$ is convex and both $\mathcal{Q}(\pi)$, $\mathcal{Q}(\rho)$ are finite. By Lemma 4.9 and (5.1), $\partial_t \pi_\varepsilon$ belong to $\mathcal{H}^{-1}(\chi(\pi_\varepsilon))$. Thus all assumptions of Lemma 4.6 are fulfilled so that $I_T(\pi_\varepsilon|\gamma) < \infty$. Since $\pi$ belongs to $\mathcal{F}_0$ and $\rho$ is the solution of the hydrodynamic equation, there exists $\delta_1 > 0$, independent of $\varepsilon$, such that $\pi_\varepsilon$ follows the hydrodynamic path on an interval $[0, \delta_1]$. Finally, by Theorem 3.3.5 in [25] and the Nash estimate, the unique solution of the hydrodynamic equation (4.2), denoted here by $\rho$, is bounded below by a strictly positive constant and above by a constant strictly smaller than 1 in any compact subset of $(0, T] \times [-1, 1]$. Hence for each $\delta_2 > 0$ there exists $a > 0$ such that $a \leq \pi_\varepsilon(t, u) \leq 1 - a$ for $\delta_2 \leq t \leq T$. This proves the claim.

Since $\pi_\varepsilon$ converges to $\pi$ as $\varepsilon \downarrow 0$, to conclude the proof of the lemma, we have to show that $\lim_{\varepsilon \to 0} I_T(\pi_\varepsilon|\gamma) = I_T(\pi|\gamma)$. To this end we verify the assumptions of Lemma 5.2. Let $\tilde{P}_\varepsilon = (1-\varepsilon)P_\pi + \varepsilon P_\rho$ and note that $\tilde{P}_\varepsilon$ is not equal to $P_{\pi_\varepsilon}$ because it does not have mean zero. To fulfill this condition, let

$$(5.2) \qquad P_\varepsilon(t, \cdot) = \tilde{P}_\varepsilon(t, \cdot) - \frac{\langle \tilde{P}_\varepsilon(t, \cdot) \chi(\pi_\varepsilon(t, \cdot))^{-1} \rangle}{\langle \chi(\pi_\varepsilon(t, \cdot))^{-1} \rangle} .$$

This expression is well defined because $\pi_\varepsilon$ is bounded away from 0 and 1 on $(0, T]$. Of course, by definition of $\tilde{P}_\varepsilon$, for every $H$ in $C_K^\infty(\Omega_T)$ $\langle\!\langle \pi_\varepsilon, \partial_t H \rangle\!\rangle = \langle\!\langle P_\varepsilon, \nabla H \rangle\!\rangle$ so that $P_\varepsilon$ is the function $P_{\pi_\varepsilon}$ defined at the beginning of this section.

Clearly, as $\varepsilon \downarrow 0$, $\pi_\varepsilon$, $\tilde{P}_\varepsilon$, $\nabla \pi_\varepsilon$ converge a.e. to $\pi$, $P$, $\nabla \pi$, respectively. We claim that $P_\varepsilon$ also converges a.e. to $P$. To prove this statement, it is enough to show that the second term on the right-hand side of (5.2) vanishes as $\varepsilon \downarrow 0$ for almost all $0 < t \leq T$. This is proved in several steps.

We first show that $\langle \chi(\pi_\varepsilon(t, \cdot))^{-1} \rangle$ converges to $\langle \chi(\pi(t, \cdot))^{-1} \rangle$ for all $0 \leq t \leq T$. Fix $0 \leq t \leq T$. On the one hand, since $\pi_\varepsilon$ converges to $\pi$ and since $\chi$ is continuous, by Fatou's lemma, $\langle \chi(\pi(t))^{-1} \rangle \leq \liminf_{\varepsilon \to 0} \langle \chi(\pi_\varepsilon(t))^{-1} \rangle$. If $\langle \chi(\pi(t))^{-1} \rangle = \infty$, we obtained the sought convergence.

Assume that $\langle \chi(\pi(t))^{-1} \rangle < \infty$. Since $\chi$ and $\chi_0$ are equivalent, this means that $\langle \chi_0(\pi(t))^{-1} \rangle < \infty$. The concavity of $\chi_0$ shows that

$$\frac{1}{\chi_0(\pi_\varepsilon)} \leq \frac{1}{(1-\varepsilon)\chi_0(\pi) + \varepsilon\chi_0(\rho)} \leq \frac{1}{(1-\varepsilon)\chi_0(\pi)}.$$

Thus if $\Xi$ stands for the bounded continuous function $\chi_0/\chi$,

$$\limsup_{\varepsilon \to 0} \Big\langle \frac{1}{\chi(\pi_\varepsilon(t))} \Big\rangle = \limsup_{\varepsilon \to 0} \Big\langle \frac{\Xi((\pi_\varepsilon(t))}{\chi_0(\pi_\varepsilon(t))} \Big\rangle \leq \limsup_{\varepsilon \to 0} \frac{1}{(1-\varepsilon)} \Big\langle \frac{\Xi((\pi_\varepsilon(t))}{\chi_0(\pi(t))} \Big\rangle.$$

Since we assumed $\chi_0(\pi(t))^{-1}$ to be integrable and since $\Xi$ is a continuous bounded function, the previous limit is equal to $\langle \Xi((\pi(t))\chi_0(\pi(t))^{-1} \rangle = \langle \chi(\pi(t))^{-1} \rangle$.

We now examine the term $\langle \tilde{P}_\varepsilon(t,\cdot)\chi(\pi_\varepsilon(t,\cdot))^{-1}\rangle$. Since $P_\rho(t,\cdot) = D(\rho_t)\nabla\rho_t - (E/2)\chi(\rho_t) + F(t)$ for some function $F\colon [0,T] \to \mathbb{R}$, $P_\rho$ is bounded for every $t > 0$. If $\langle \chi(\pi(t,\cdot))^{-1}\rangle$ is finite, an argument similar to the one presented in the previous paragraphs shows that $\langle \tilde{P}_\varepsilon(t,\cdot)\chi(\pi_\varepsilon(t,\cdot))^{-1}\rangle$ converges to $\langle P_\pi(t,\cdot)\chi(\pi(t,\cdot))^{-1}\rangle$. This expression vanishes by definition of $P_\pi$.

Suppose now that $\langle \chi(\pi(t,\cdot))^{-1}\rangle = \infty$. By the Schwarz inequality,

$$\left(\frac{\langle P_\pi(t,\cdot)\chi(\pi_\varepsilon(t,\cdot))^{-1}\rangle}{\langle \chi(\pi_\varepsilon(t,\cdot))^{-1}\rangle}\right)^2 \le \frac{\langle P_\pi(t,\cdot)^2\chi(\pi_\varepsilon(t,\cdot))^{-1}\rangle}{\langle \chi(\pi_\varepsilon(t,\cdot))^{-1}\rangle}.$$

We have already seen that the denominator diverges while the numerator remains bounded by $C_0\langle P_\pi(t,\cdot)^2\chi(\pi(t,\cdot))^{-1}\rangle$ for some finite constant $C_0$. This expression is finite a.s. in $t$ because $P_\pi$ belongs to $\mathbb{L}^2(\chi(\pi)^{-1})$. On the other hand, $\langle P_\rho(t,\cdot)\chi(\pi_\varepsilon(t,\cdot))^{-1}\rangle$ is bounded by $C_0\langle \chi(\pi_\varepsilon(t,\cdot))^{-1}\rangle$ for some finite constant $C_0$ which depends on $\rho$. Putting together the previous two assertions, we obtain that the second term on the right-hand side of (5.2) vanishes if $\langle \chi(\pi(t,\cdot))^{-1}\rangle = \infty$. This concludes the proof that $P_\varepsilon$ converges a.e. to $P$.

To prove hypothesis (2) in Lemma 5.2, first note that it is enough to show that $\{\nabla\pi_\varepsilon\}^2/\chi_0(\pi_\varepsilon)$, $P_{\pi_\varepsilon}^2/\chi_0(\pi_\varepsilon)$ are uniformly integrable sequences. By (5.1) and Lemma 5.3, there exists a convex increasing function $\Psi$ such that

$$\int_0^T dt \int_{-1}^1 du\, \Psi\left(\frac{\{\nabla\pi\}^2}{\chi_0(\pi)}\right) < \infty, \qquad \int_0^T dt \int_{-1}^1 du\, \Psi\left(\frac{\{\nabla\rho\}^2}{\chi_0(\rho)}\right) < \infty.$$

By the Schwarz inequality, $\{\nabla\pi_\varepsilon(t,u)\}^2$, which is equal to the right-hand side of the next equation, is bounded above by

$$\left\{(1-\varepsilon)\frac{\nabla\pi(s,u)}{\sqrt{\chi_0(\pi(s,u))}}\sqrt{\chi_0(\pi(s,u))} + \varepsilon\frac{\nabla\rho(s,u)}{\sqrt{\chi_0(\rho(s,u))}}\sqrt{\chi_0(\rho(s,u))}\right\}^2$$

$$\le \{(1-\varepsilon)\chi_0(\pi(s,u)) + \varepsilon\chi_0(\rho(s,u))\}$$

$$\times \left\{(1-\varepsilon)\frac{\{\nabla\pi(s,u)\}^2}{\chi_0(\pi(s,u))} + \varepsilon\frac{\{\nabla\rho(s,u)\}^2}{\chi_0(\rho(s,u))}\right\}.$$

By the concavity of $\chi_0$ and Jensen's inequality, this expression is less than or equal to

$$\chi_0(\pi_\varepsilon(s,u))\left\{(1-\varepsilon)\frac{\{\nabla\pi(s,u)\}^2}{\chi_0(\pi(s,u))} + \varepsilon\frac{\{\nabla\rho(s,u)\}^2}{\chi_0(\rho(s,u))}\right\}.$$

Hence since $\Psi$ is increasing and convex,

$$\int_0^T dt \int_{-1}^1 du\, \Psi\left(\frac{\{\nabla\pi_\varepsilon\}^2}{\chi_0(\pi_\varepsilon)}\right) \le (1-\varepsilon)\int_0^T dt \int_{-1}^1 du\, \Psi\left(\frac{\{\nabla\pi\}^2}{\chi_0(\pi)}\right)$$

$$+ \varepsilon \int_0^T dt \int_{-1}^1 du\, \Psi\left(\frac{\{\nabla\rho\}^2}{\chi_0(\rho)}\right).$$

Thus by Lemma 5.3, the sequence $\{\nabla \pi_\varepsilon\}^2/\chi_0(\pi_\varepsilon)$ is uniformly integrable. We may proceed in a similar way to prove the uniform integrability of $P_{\pi_\varepsilon}^2/\chi_0(\pi_\varepsilon)$.

To prove assumption (3), note that it is enough to show that $\lim_{\varepsilon \to 0}\langle \chi(\pi_\varepsilon(t))^{-1}\rangle = \langle \chi(\pi(t))^{-1}\rangle$ for almost all $0 \leq t \leq T$ and to apply the dominated convergence theorem because $\chi$ is bounded, but this has already been proved. $\square$

Denote by $\mathcal{F}_2$ the set of trajectories $\pi$ in $\mathcal{F}_1$ for which there exists $\delta_1$, $\delta_2 > 0$ such that: $\pi$ follow the hydrodynamic path in the time interval $[0, \delta_1]$, is constant in time in the interval $[\delta_1, \delta_1 + \delta_2]$, and $\pi$ is smooth in time in the time interval $(\delta_1, T]$.

LEMMA 5.6. *The set $\mathcal{F}_2$ is $I_T(\cdot|\gamma)$-dense.*

PROOF. Fix a trajectory $\pi$ in $\mathcal{F}_1$. Assume that $\pi$ follows the hydrodynamic equation in the time interval $[0, 2a]$ for some $a > 0$. Let $a \leq t_0 \leq 2a$ be such that $\langle (\nabla \pi_{t_0})^2/\chi(\pi_{t_0})\rangle < \infty$. This is possible because $\pi$ follows the hydrodynamic path $\rho$ in the time interval $[0, 2a]$ and $\mathcal{Q}(\rho) < \infty$. For $0 < \varepsilon < a$, let $T - 2\varepsilon < T_\varepsilon < T - \varepsilon$ such that

$$\left\langle \frac{(\nabla \pi_{T_\varepsilon})^2}{\chi(\pi_{T_\varepsilon})} \right\rangle \leq \frac{1}{\varepsilon}\int_{T-2\varepsilon}^{T-\varepsilon} \left\langle \frac{(\nabla \pi_t)^2}{\chi(\pi_t)} \right\rangle dt + \varepsilon,$$

$\pi(T_\varepsilon, \pm 1) = \rho_\pm$. Roughly speaking, the profile $\pi_{T_\varepsilon}$ minimizes locally the energy and has the correct boundary conditions. The latter condition can be achieved because $\pi$ belongs to $D_\gamma$ and $\mathcal{Q}(\pi) < \infty$ so that $\pi(t, \cdot)$ is continuous and $\pi(t, \pm 1) = \rho_\pm$ for almost all $t$. Clearly, $\langle (\nabla \pi_{T_\varepsilon})^2/\chi(\pi_{T_\varepsilon})\rangle \leq \varepsilon^{-1}\mathcal{Q}(\pi) + \varepsilon$.

Define the path $\tilde{\pi}_\varepsilon$ by

$$\tilde{\pi}_\varepsilon(t, \cdot) = \begin{cases} \pi(t, \cdot), & \text{for } 0 \leq t \leq t_0, \\ \pi(t_0, \cdot), & \text{for } t_0 \leq t \leq t_0 + \varepsilon, \\ \pi(t - \varepsilon, \cdot), & \text{for } t_0 + \varepsilon \leq t \leq T_\varepsilon + \varepsilon, \\ \pi(T_\varepsilon, \cdot), & \text{for } T_\varepsilon + \varepsilon \leq t \leq T. \end{cases}$$

We claim that $\tilde{\pi}_\varepsilon$ belongs to $\mathcal{F}_1$ for each $\varepsilon > 0$. Moreover, as $\varepsilon \downarrow 0$, $\tilde{\pi}_\varepsilon$, $I_T(\tilde{\pi}_\varepsilon|\gamma)$ converge to $\pi$, $I_T(\pi|\gamma)$, respectively.

By construction $\tilde{\pi}_\varepsilon$ belongs to $D_\gamma$ and follows the hydrodynamic equation in the time interval $[0, t_0]$. It is also bounded below by a strictly positive constant and above by a constant strictly less than 1 on each time interval $[\delta, T]$, $\delta > 0$ because $\pi$ is as well. Moreover, $\mathcal{Q}(\tilde{\pi}_\varepsilon)$ is uniformly bounded by $\mathcal{Q}(\pi) + \varepsilon\langle(\nabla \pi_{t_0})^2/\chi(\pi_{t_0})\rangle + \varepsilon\langle(\nabla \pi_{T_\varepsilon})^2/\chi(\pi_{T_\varepsilon})\rangle \leq 2\mathcal{Q}(\pi) + O(\varepsilon)$. Finally, to show that $\hat{I}_T(\tilde{\pi}_\varepsilon|\gamma)$ is uniformly bounded, for $0 \leq t_0 < t_1 \leq T$, denote by $\hat{I}_{[t_0, t_1]}(\cdot|\gamma)$ the contribution of the time interval $[t_0, t_1]$ to the rate function $\hat{I}_T(\tilde{\pi}_\varepsilon|\gamma)$. On the one hand, $\hat{I}_{[0,t_0]\cup[t_0+\varepsilon,T_\varepsilon+\varepsilon]}(\tilde{\pi}_\varepsilon|\gamma) \le \hat{I}_T(\pi|\gamma)$. On the other hand, $\hat{I}_{[t_0,t_0+\varepsilon]}(\tilde{\pi}_\varepsilon|\gamma) \le C_0\varepsilon\langle(\nabla\pi_{t_0})^2/\chi(\pi_{t_0})\rangle$ and $\hat{I}_{[T_\varepsilon+\varepsilon,T]}(\tilde{\pi}_\varepsilon|\gamma) \le C_0\varepsilon\langle(\nabla\pi_{T_\varepsilon})^2/\chi(\pi_{T_\varepsilon})\rangle$ for some finite constant $C_0$ depending only on $E$, $D(\cdot)$, $\chi(\cdot)$. By definition of $T_\varepsilon$, this expression is bounded by

$$C_0 \int_{T-2\varepsilon}^{T-\varepsilon} \left\langle \frac{(\nabla\pi_t)^2}{\chi(\pi_t)} \right\rangle dt + C_0\varepsilon^2,$$

which vanishes as $\varepsilon \downarrow 0$ because $\pi$ has finite energy. This proves that $\tilde{\pi}_\varepsilon$ belongs to $\mathcal{F}_1$ for each $\varepsilon > 0$.

Since $\pi$ belongs to $C([0,T],\mathcal{M})$, $\tilde{\pi}_\varepsilon$ converges to $\pi$ as $\varepsilon \downarrow 0$. By lower semicontinuity, to prove that $\hat{I}_T(\tilde{\pi}_\varepsilon|\gamma)$ converges to $\hat{I}_T(\pi|\gamma)$ it is enough to show that $\limsup_{\varepsilon\to 0} \hat{I}_T(\tilde{\pi}_\varepsilon|\gamma) \le \hat{I}_T(\pi|\gamma)$. This follows from the bound on $\hat{I}_T(\tilde{\pi}_\varepsilon|\gamma)$ obtained in the previous paragraph.

In conclusion, to prove the lemma, it is enough to show that every path $\pi$ in $\mathcal{F}_1$ which follows the hydrodynamic equation in a time interval $[0,b]$, is constant in the time intervals $[b,b+a]$, $[T-a,T]$, can be approximated by a sequence in $\mathcal{F}_2$. Fix such a path $\pi$.

Denote by $\iota$ a smooth, positive function with support contained in $[0,1]$ and integral equal to 1. For $\varepsilon > 0$, let $\iota_\varepsilon(t) = \varepsilon^{-1}\iota(t\varepsilon^{-1})$ be a smooth approximation of the identity. For $\varepsilon < a$, consider the path $\pi_\varepsilon$ defined by

$$\pi_\varepsilon(t,\cdot) = \begin{cases} \pi(t,\cdot), & \text{for } 0 \le t \le b, \\ \int_\mathbb{R} ds\,\iota_\varepsilon(s)\pi(t+s,\cdot), & \text{for } b \le t \le T, \end{cases}$$

where we extended the definition of $\pi$ to the time interval $[T,\infty)$ by setting $\pi(t,u) = \pi(T,u)$ for $t \ge T$.

We claim that $\pi_\varepsilon$ belongs to $\mathcal{F}_2$ for each $\varepsilon < a$. By construction it belongs to $D_\gamma$, follows the hydrodynamic equation in the time interval $[0,b]$ and is bounded below by a strictly positive constant and above by a constant strictly less than 1 on each time interval $[\delta',T]$, $\delta' > 0$. It is constant in time in the interval $[b, b+a-\varepsilon]$ and smooth in time in the interval $(b,T]$ because we chose a smooth approximation $\iota_\varepsilon$ of the identity. By convexity of $\mathcal{Q}$, $\mathcal{Q}(\pi_\varepsilon) \le \mathcal{Q}(\pi) + O(\varepsilon)$ because $\mathcal{Q}(\pi_b^*)$ and $\mathcal{Q}(\pi_T^*)$ are finite. Here $\pi_b^*$, $\pi_T^*$ are the trajectories constant in time and equal at each time to $\pi_b$ and $\pi_T$, respectively.

It remains to show that $\pi_\varepsilon$, $\hat{I}_T(\pi_\varepsilon|\gamma)$ converge to $\pi$, $\hat{I}_T(\pi|\gamma)$, respectively. We rely on Lemma 5.2.

Recall from (4.12) that there exists $P = P_\pi$ in $\mathbb{L}^2(\chi(\pi)^{-1})$ such that $\langle\!\langle\pi, \partial_t H\rangle\!\rangle = \langle\!\langle P, \nabla H\rangle\!\rangle$ for every $H$ in $C_K^\infty(\Omega_T)$. An elementary computation shows that $\langle\!\langle\pi_\varepsilon, \partial_t H\rangle\!\rangle = \langle\!\langle\tilde{P}_\varepsilon, \nabla H\rangle\!\rangle$ if we define $\tilde{P}_\varepsilon$ by

$$\tilde{P}_\varepsilon(t,\cdot) = \begin{cases} P(t,\cdot), & \text{for } 0 \le t \le b, \\ \int_\mathbb{R} ds\,\iota_\varepsilon(s)P(t+s,u), & \text{for } b \le t \le T. \end{cases}$$

Note that $\tilde{P}_\varepsilon$ belongs to $\mathbb{L}^2(\chi(\pi_\varepsilon)^{-1})$ because $\tilde{P}_\varepsilon$, $\pi_\varepsilon$ coincide with $P$, $\pi$ on the interval $[0, b]$ because $\pi$ is bounded below by a strictly positive constant and above by a constant strictly less than 1 on $[b, T]$ so that the denominator $\chi^{-1}$ is irrelevant and because the square function is convex. In particular, by Lemma 4.6, $\hat{I}_T(\pi_\varepsilon | \gamma)$ is finite for every $\varepsilon > 0$.

Let
$$P_\varepsilon(t, \cdot) = \tilde{P}_\varepsilon(t, \cdot) - \frac{\langle \tilde{P}_\varepsilon(t, \cdot)\chi(\pi_\varepsilon(t, \cdot))^{-1}\rangle}{\langle \chi(\pi_\varepsilon(t, \cdot))^{-1}\rangle} \mathbf{1}\{\langle \chi(\pi_\varepsilon(t, \cdot))^{-1}\rangle < \infty\}$$
for $P_\varepsilon$ to have mean zero. We may remove the indicator function because $\pi_\varepsilon$ is bounded below and above on $(0, T]$.

Since $\pi$, $\nabla\pi$, $P$ belong to $L^1(\Omega_T)$, $\pi_\varepsilon$, $\nabla\pi_\varepsilon$, $\tilde{P}_\varepsilon$ converge in $L^1(\Omega_T)$ to $\pi$, $\nabla\pi$, $P$, respectively. Taking subsequences, if necessary, we obtain a.e. convergence and convergence of $\tilde{P}_\varepsilon(t, \cdot)$ to $P(t, \cdot)$ in $L^1(\Omega)$ for a.e. $0 \leq t \leq T$. Since $\pi_\varepsilon$ is bounded away from 0, 1, it is not difficult to show that $P_\varepsilon$ also converges a.e. to $P$.

Since $D$ is bounded and $\chi$ is equivalent to $\chi_0$, to prove the uniform integrability of $\{D(\pi_\varepsilon)\nabla\pi_\varepsilon\}^2/\chi(\pi_\varepsilon)$, it is enough to estimate $\{\nabla\pi_\varepsilon\}^2/\chi_0(\pi_\varepsilon)$. Recall that $\mathcal{Q}(\pi_T^*) < \infty$. By Lemma 5.3, there exists a convex increasing function $\Psi$ such that
$$(5.3) \quad \int_0^T dt \int_{-1}^1 du\, \Psi\left(\frac{\{\nabla\pi\}^2}{\chi_0(\pi)}\right) < \infty, \qquad \int_{-1}^1 du\, \Psi\left(\frac{\{\nabla\pi_T\}^2}{\chi_0(\pi_T)}\right) < \infty.$$

By the Schwarz inequality, the concavity of $\chi_0$ and the Jensen inequality, for $t \geq b$,
$$\{\nabla\pi_\varepsilon(t, u)\}^2 = \left(\int ds\, \iota_\varepsilon(s-t) \frac{\nabla\pi(s, u)}{\chi_0(\pi(s,u))^{1/2}} \chi_0(\pi(s,u))^{1/2}\right)^2$$
$$\leq \chi_0(\pi_\varepsilon(t, u)) \int ds\, \iota_\varepsilon(s-t) \frac{\{\nabla\pi(s, u)\}^2}{\chi_0(\pi(s, u))}.$$

Hence, by the Jensen inequality, since $\Psi$ is convex and increasing,
$$\int_b^T dt \int_{-1}^1 du\, \Psi\left(\frac{\{\nabla\pi_\varepsilon\}^2}{\chi_0(\pi_\varepsilon)}\right) \leq \int_b^T dt \int_{-1}^1 du \int ds\, \iota_\varepsilon(s-t) \Psi\left(\frac{\{\nabla\pi(s, u)\}^2}{\chi_0(\pi(s, u))}\right).$$

Integrating in $t$ the right-hand side, we get that
$$\int_b^T dt \int_{-1}^1 du\, \Psi\left(\frac{\{\nabla\pi_\varepsilon\}^2}{\chi_0(\pi_\varepsilon)}\right)$$
$$\leq \int_b^T ds \int_{-1}^1 du\, \Psi\left(\frac{\{\nabla\pi\}^2}{\chi_0(\pi)}\right) + \varepsilon \int_{-1}^1 du\, \Psi\left(\frac{\{\nabla\pi_T\}^2}{\chi_0(\pi_T)}\right) < \infty.$$

It remains to add the piece corresponding to the time interval $[0, b]$ to obtain, in view of (5.3) and Lemma 5.3, the uniform integrability of $\{\nabla\pi_\varepsilon\}^2/\chi_0(\pi_\varepsilon)$

and, therefore, the one of $\{D(\pi_\varepsilon)\nabla\pi_\varepsilon\}^2/\chi(\pi_\varepsilon)$. The same argument shows the uniform integrability of $P_\varepsilon^2/\chi(\pi_\varepsilon)$.

To prove the third condition of Lemma 5.2, it is enough to show that $\lim_{\varepsilon\to 0}\langle\chi(\pi_\varepsilon(t))^{-1}\rangle = \langle\chi(\pi(t))^{-1}\rangle$ for all $0 < t \le T$. Fix $t > 0$. Since $\pi$ belongs to $\mathcal{F}_1$, there exists $\delta > 0$ and $\varepsilon_0 > 0$ such that $\delta \le \pi_\varepsilon(t) \le 1-\delta$ for all $\varepsilon < \varepsilon_0$. Since $\pi_\varepsilon(t)$ converges to $\pi(t)$ a.s. as $\varepsilon \downarrow 0$, condition (3) follows from the dominated convergence theorem. □

Note that each path $\pi_\varepsilon$ defined in the proof of this lemma is continuous on $(0,T]\times[-1,1]$. The path $\pi_\varepsilon$ is continuous on $(0,b]\times[-1,1]$ because it follows the hydrodynamic equation. The continuity can be extended to $[b,b+a-\varepsilon]\times[-1,1]$ because $\pi_\varepsilon(t,\cdot)$ is constant and equal to $\pi(b,\cdot)$ in this interval. By construction, $\pi_\varepsilon$ is continuous in time on $(0,T]\times[-1,1]$. On the other hand, if we denote by $*$ the convolution, since for $t \ge b$ $\pi_\varepsilon(t,\cdot) = (\pi * \iota_\varepsilon)(t,\cdot)$, by the Schwarz inequality, for $b \le t \le T$,

$$\langle(\nabla\pi_\varepsilon(t,\cdot))^2\rangle \le \int_0^\varepsilon \iota_\varepsilon(s)\langle(\nabla\pi_{t+s})^2\rangle\,ds \le \frac{C_1}{\varepsilon}\int_t^{t+\varepsilon}\langle(\nabla\pi_s)^2\rangle\,ds \le \frac{2C_1}{\varepsilon}\mathcal{Q}(\pi).$$

In the last step we used the fact that $\pi$ is constant in the interval $[T-a,\infty)$ to estimate the piece of the integral in the interval $[T,T+\varepsilon]$ by the energy of $\pi$. It follows from this energy estimate that $\pi_\varepsilon$ is uniformly continuous in space in $[b,T]\times[-1,1]$. In particular, we could have required paths in $\mathcal{F}_2$ to be continuous, but we do not need this property in the sequel.

Let $\mathcal{F}_3$ be the set of trajectories $\pi$ in $\mathcal{F}_1$ (not $\mathcal{F}_2$) for which there exists $\delta > 0$ such that: $\pi$ follows the hydrodynamic path in the time interval $[0,\delta]$, is continuous on $(0,T]\times[-1,1]$ and smooth on $[\delta,T]\times[-1,1]$. Note that $\mathcal{F}_3$ corresponds to the set $\mathcal{D}_{T,\gamma}^\circ$ introduced in Definition 3.6. For a path $\pi$ in $\mathcal{F}_3$, denote by $\mathfrak{t}(\pi)$ the positive time at which smoothness may be violated. In the previous description, $\mathfrak{t}(\pi) = b$ once one recalls that $\pi$ is also smooth in $(0,\delta)\times[-1,1]$ because it follows the hydrodynamic equation.

PROOF OF THEOREM 5.1. Fix a trajectory $\pi$ in $\mathcal{F}_2$. In view of the previous lemma, it is enough to show that there exists a sequence $\{\pi_\varepsilon\}$ in $\mathcal{F}_3$ such that $\pi_\varepsilon$, $I_T(\pi_\varepsilon|\gamma)$ converge to $\pi$, $I_T(\pi|\gamma)$, respectively.

For $\varepsilon > 0$, denote by $R_\varepsilon^D$, $R_\varepsilon^N : [-1,1]^2 \to \mathbb{R}_+$ the resolvent of the Dirichlet, respectively, Neumann Laplacian: $R_\varepsilon^{D/N} = (I - \varepsilon\Delta_{D/N})^{-1}$, where $\Delta_D$, $\Delta_N$ stand for the Laplacian with Dirichlet, Neumann boundary conditions. An elementary computation gives an explicit form for the resolvent transcribed below:

$$R_\varepsilon^D(u,v) = \frac{\sqrt{\lambda}}{\sinh\{2\sqrt{\lambda}\}}\begin{cases}\sinh\{\sqrt{\lambda}(1+u)\}\sinh\{\sqrt{\lambda}(1-v)\}, & \text{if } u \le v, \\ \sinh\{\sqrt{\lambda}(1-u)\}\sinh\{\sqrt{\lambda}(1+v)\}, & \text{otherwise;}\end{cases}$$

$$R_\varepsilon^N(u,v) = \frac{\sqrt{\lambda}}{\sinh\{2\sqrt{\lambda}\}}\begin{cases}\cosh\{\sqrt{\lambda}(1+u)\}\cosh\{\sqrt{\lambda}(1-v)\}, & \text{if } u \le v, \\ \cosh\{\sqrt{\lambda}(1-u)\}\cosh\{\sqrt{\lambda}(1+v)\}, & \text{otherwise;}\end{cases}$$

where $\lambda = \varepsilon^{-1}$. In contrast with $R^D$, $R^N$ is a probability kernel.

Since $\pi$ belongs to $\mathcal{F}_2$, there exists $0 < a < b \leq T$, such that $\pi$ follows the hydrodynamic equation in the time interval $[0, a]$ and is constant in the time interval $[a, b]$. Let $j : \mathbb{R} \to [0, 1]$ be a smooth nondecreasing function such that $j(t) = 0, 1$, for $t \leq 0, t \geq 1$, respectively. For $0 < \varepsilon < b - a$, let $j_\varepsilon(t) = \varepsilon j(t\varepsilon^{-1})$, $\beta_\varepsilon(t) = j_\varepsilon(t - a)$ and

$$\pi_\varepsilon(t, \cdot) = \begin{cases} \pi(t, \cdot), & \text{for } 0 \leq t \leq a, \\ \rho^* + R^D_{\beta_\varepsilon(t)}(\pi(t, \cdot) - \rho^*), & \text{for } a < t \leq T, \end{cases}$$

where $\rho^*$ is the linear profile $\rho_-(1-u)/2 + \rho_+(1+u)/2$.

We claim that $\pi_\varepsilon$ belongs to $\mathcal{F}_3$. Since $\pi$ belongs to $\mathcal{F}_2$, by construction, $\pi_\varepsilon$ belongs to $D_\gamma$, $\pi_\varepsilon$ follows the hydrodynamic equation on the time interval $[0, a]$ and $\pi_\varepsilon$ is smooth in space and time on the interval $(a, T]$ and continuous on $[a, T] \times [-1, 1]$. We prove at the end of the lemma that there exists a finite constant $C_0$ such that

(5.4) $$\frac{\chi_0(R^N_{\beta_\varepsilon(t)} \pi(t, u))}{\chi_0(\pi_\varepsilon(t, u))} \leq C_0$$

for all $t > a$, $-1 \leq u \leq 1$, $\varepsilon > 0$. Since $R^N_\varepsilon$ is the resolvent of the Laplacian with reflecting boundary conditions, $\inf_{u \in [-1,1]} f(u) \leq (R^N_\varepsilon f)(v) \leq \sup_{u \in [-1,1]} f(u)$ for every $v$ in $[-1, 1]$, $\varepsilon > 0$. In particular, there exists $\delta > 0$ such that $\delta \leq (R^N_\varepsilon \pi)(t, u) \leq 1 - \delta$ for $t \geq a$ because $\pi$ belongs to $\mathcal{F}_2$. It follows from (5.4) that the same holds for $\pi_\varepsilon$. Therefore, $\pi_\varepsilon$ is bounded away from 0 and 1 in the time interval $(a, T]$. Since $\pi_\varepsilon$ follows the hydrodynamic equation in $[0, a]$, for every $\delta > 0$, there exists $\varphi > 0$, independent of $\varepsilon$, such that $\varphi \leq \pi_\varepsilon \leq 1 - \varphi$ in the time interval $[\delta, T]$ for all $\varepsilon > 0$.

We examine in this paragraph the energy of $\pi_\varepsilon$. An elementary computation shows that $\nabla \pi_\varepsilon = R^N_{\beta_\varepsilon} \nabla \pi$. Therefore, by the Schwarz inequality, for $t \geq a$,

$$\langle (\nabla \pi_\varepsilon)^2 \rangle = \langle (R^N_{\beta_\varepsilon} \nabla \pi)^2 \rangle \leq \langle R^N_{\beta_\varepsilon} (\nabla \pi)^2 \rangle = \langle (\nabla \pi)^2 \rangle.$$

Since $\pi_\varepsilon = \pi$ on the time interval $[0, a]$ and since both paths are bounded away from 0 and 1 in the time interval $[a, T]$, there exists a finite constant $C_0 = C_0(\pi)$ such that $\mathcal{Q}(\pi_\varepsilon) \leq C_0 \mathcal{Q}(\pi)$ uniformly over $\varepsilon > 0$.

To show that $\hat{I}_T(\pi_\varepsilon | \gamma)$ is bounded, recall from (4.12) that there exists $P = P_\pi$ in $\mathbb{L}^2(\chi(\pi)^{-1})$ such that $\langle\!\langle \pi, \partial_t H \rangle\!\rangle = \langle\!\langle P, \nabla H \rangle\!\rangle$ for every $H$ in $C_K^\infty(\Omega_T)$. Fix such a function $H$. A straightforward computation using the relation $\nabla R^D_\varepsilon = R^N_\varepsilon \nabla$ shows that

(5.5) $$\langle\!\langle \pi_\varepsilon, \partial_t H \rangle\!\rangle = \langle\!\langle R^N_{\beta_\varepsilon(t)} P, \nabla H \rangle\!\rangle - \langle\!\langle \pi - \rho^*, (\partial_t R^D_{\beta_\varepsilon(t)}) H \rangle\!\rangle.$$

In the first term on the right-hand side, it must be understood that $R^N_{\beta_\varepsilon(t)} P = P$ for $t \leq a$ and in the second term that the time derivative concerns $R^D_{\beta_\varepsilon(t)}$ exclusively.

We claim that the second term on the right-hand side of (5.5) is negligible in the sense that the linear functional $\ell_\varepsilon : C_K^\infty(\Omega_T) \to \mathbb{R}$, $\ell_\varepsilon(H) = \langle\!\langle \pi - \rho^*, (\partial_t R^D_{\beta_\varepsilon(t)}) H \rangle\!\rangle$, is bounded in $\mathcal{H}^{-1}(\chi(\pi_\varepsilon))$ by a constant which vanishes as $\varepsilon \downarrow 0$. More precisely, there exists a constant $C_\varepsilon(\pi)$, which vanishes as $\varepsilon \downarrow 0$, such that

(5.6) $$\langle\!\langle \pi - \rho^*, (\partial_t R^D_{\beta_\varepsilon(t)}) H \rangle\!\rangle^2 \leq C_\varepsilon(\pi) \|H\|^2_{1,\chi(\pi_\varepsilon)}$$

for every $H$ in $C_K^\infty(\Omega_T)$. Indeed, since $\partial_t R^D_{\beta_\varepsilon(t)} = \beta'_\varepsilon(t) \Delta_D (R^D_{\beta_\varepsilon(t)})^2$ where $R^2$ stands for the composition of the operator $R$ with itself, by the Schwarz inequality, the left-hand side of (5.6) is bounded above by

$$\int_0^T [\beta'_\varepsilon(t)]^2 \langle \{\nabla[\pi_t - \rho^*]\}^2 \rangle \, dt \int_a^{a+\varepsilon} \langle (\nabla H_t)^2 \rangle \, dt.$$

Since $\pi_\varepsilon$ is bounded away from 0 and 1 in the time interval $[a, a+\varepsilon]$, we may include a factor $\chi(\pi_\varepsilon(t, \cdot))$ in the second integral paying the price of a constant $C_0 = C_0(\pi)$ and then extend the time integral to the interval $[0, T]$. Since $\beta'_\varepsilon(t)$ vanishes outside the interval $[a, a+\varepsilon]$ and is bounded uniformly in time and in $\varepsilon > 0$, the previous expression is less than or equal to

$$C_0(\pi) \int_a^{a+\varepsilon} \langle \{\nabla[\pi_t - \rho^*]\}^2 \rangle \, dt \int_0^T \langle (\nabla H_t)^2 \chi(\pi_\varepsilon(t, \cdot)) \rangle \, dt.$$

Since $\pi$ and $\rho^*$ have finite energy, the first integral vanishes as $\varepsilon \downarrow 0$. This proves (5.6).

We claim that $\partial_t \pi_\varepsilon$ belongs to $\mathcal{H}^{-1}(\chi(\pi_\varepsilon))$. In view of (5.5), (5.6), it is enough to show that $R^N_{\beta_\varepsilon(t)} P$ belongs to $\mathbb{L}^2(\chi(\pi_\varepsilon)^{-1})$. This follows from the identities $\pi_\varepsilon = \pi$, $R^N_{\beta_\varepsilon(t)} P = P$ for $t \leq a$, and the observation that $\langle (R^N_{\beta_\varepsilon(t)} P)^2 \rangle \leq \langle P^2 \rangle$, together with the fact that $\pi_\varepsilon$, $\pi$ are bounded away from 0, 1 in the time interval $[a, T]$, uniformly in $\varepsilon > 0$.

Since $\partial_t \pi_\varepsilon$ belongs to $\mathcal{H}^{-1}(\chi(\pi_\varepsilon))$ and $\pi_\varepsilon$ has finite energy, by Lemma 4.6, $\hat{I}_T(\pi_\varepsilon | \gamma)$ is finite. To conclude the proof of the lemma, it remains to show that $\pi_\varepsilon$, $\hat{I}_T(\pi_\varepsilon | \gamma)$ converge to $\pi$, $\hat{I}_T(\pi | \gamma)$.

By construction, $\pi_\varepsilon$ converges to $\pi$. Let

$$P_\varepsilon = R^N_{\beta_\varepsilon} P - \frac{\langle R^N_{\beta_\varepsilon} P \chi(\pi_\varepsilon)^{-1} \rangle}{\langle \chi(\pi_\varepsilon)^{-1} \rangle}.$$

By (5.5), for every $H$ in $C_K^\infty(\Omega_T)$, $\langle\!\langle \pi_\varepsilon, \partial_t H \rangle\!\rangle = \langle\!\langle P_\varepsilon, \nabla H \rangle\!\rangle + \ell_\varepsilon(H)$, where, by (5.6), $\ell_\varepsilon(H)$ is a negligible term. In particular, by Lemma 5.2, to show that $\hat{I}_T(\pi_\varepsilon | \gamma)$ converges to $\hat{I}_T(\pi | \gamma)$ we just need to check the assumptions (1), (2), (3) for $\pi_\varepsilon$, $\nabla \pi_\varepsilon$ and $P_\varepsilon$.

By definition, for $t > a$, $\pi_\varepsilon = \rho^* + R^D_{\beta_\varepsilon}(\pi - \rho^*)$, $\nabla \pi_\varepsilon = R^N_{\beta_\varepsilon} \nabla \pi$ and $R^N_{\beta_\varepsilon} P$ converge a.e. to $\pi$, $\nabla \pi$, $P$, respectively. Repeating the arguments presented in the proof of the previous lemma, we can deduce from the a.s. convergence of $R^N_{\beta_\varepsilon} P$ that $P_\varepsilon$ converges a.e. to $P$.

We show that $(\nabla \pi_\varepsilon)^2/\chi(\pi_\varepsilon)$ is uniformly integrable. Since $R^N$ is a probability kernel and since $\chi_0$ is concave, by the Schwarz and Jensen inequalities,

$$\{\nabla \pi_\varepsilon(t,u)\}^2 = \left(\int dv\, R^N_{\beta_\varepsilon}(u,v) \frac{\nabla \pi(t,v)}{\chi_0(\pi(t,v))^{1/2}} \chi_0(\pi(t,v))^{1/2}\right)^2$$

$$\leq \chi_0(R^N_{\beta_\varepsilon}\pi(t,u)) \int dv\, R^N_{\beta_\varepsilon}(u,v) \frac{\{\nabla \pi(t,v)\}^2}{\chi_0(\pi(t,v))}.$$

By (5.4) and proceeding as in the proof of the previous lemma, we end the proof of the uniformly integrability of $\{\nabla \pi_\varepsilon\}^2/\chi_0(\pi_\varepsilon)$. The same argument applies to $P_\varepsilon$.

Assumption (3) of Lemma 5.2 is proved as in the previous lemma. This concludes the proof, modulo (5.4), which we now examine. Recall that $\lambda = \varepsilon^{-1}$. Let us first consider the case $u \in [-1, -1 + \sqrt{\varepsilon}]$. By an explicit computation,

$$\pi_\varepsilon(t,u) \geq \rho^*(u) - R^D_\varepsilon \rho^*(u) = \rho_- \frac{\sinh[\sqrt{\lambda}(1-u)]}{\sinh 2\sqrt{\lambda}} + \rho_+ \frac{\sinh[\sqrt{\lambda}(1+u)]}{\sinh 2\sqrt{\lambda}}.$$

Since $\pi \leq 1$, $R^N_\varepsilon \pi \leq 1$ so that

$$\sup_{\varepsilon \leq 1} \sup_{\substack{t \in [a,T] \\ u \in [-1,-1+\sqrt{\varepsilon}]}} \frac{R^N_{\beta_\varepsilon(t)} \pi(t,u)}{\pi_\varepsilon(t,u)} \leq C_0$$

for some finite constant $C_0$. By analogous computations

$$1 - \pi_\varepsilon(t,u) \geq 1 - \rho^*(u) - R^D_\varepsilon(1-\rho^*)(u)$$

$$= (1-\rho_-) \frac{\sinh[\sqrt{\lambda}(1-u)]}{\sinh 2\sqrt{\lambda}} + (1-\rho_+) \frac{\sinh[\sqrt{\lambda}(1+u)]}{\sinh 2\sqrt{\lambda}}.$$

Since $1 - R^N_\varepsilon \pi = R^N_\varepsilon(1-\pi) \leq 1$, we get

$$\sup_{\varepsilon \leq 1} \sup_{\substack{t \in [a,T] \\ u \in [-1,-1+\sqrt{\varepsilon}]}} \frac{1 - R^N_{\beta_\varepsilon(t)} \pi(t,u)}{1 - \pi_\varepsilon(t,u)} \leq C_0$$

for some finite constant $C_0$, which yields the bound (5.4) for $u \in [-1, -1 + \sqrt{\varepsilon}]$. Of course, the same argument applies for $u \in [1 - \sqrt{\varepsilon}, 1]$.

To analyze the case $u \in A^\varepsilon = [-1 + \sqrt{\varepsilon}, 1 - \sqrt{\varepsilon}]$, we first show that there exists a constant $C_1$ such that

$$
\begin{aligned}
& \int_{-1}^{u} dv \cosh[\sqrt{\lambda}(1+v)]\pi(t,v) \\
& \qquad \leq C_1 \int_{-1}^{u} dv \sinh[\sqrt{\lambda}(1+v)]\pi(t,v)
\end{aligned}
\tag{5.7}
$$

for any $0 < \varepsilon \leq 1$, $a \leq t \leq T$, $u$ in $A^\varepsilon$.

Since $\pi$ is bounded away from 0 and 1 in the time interval $[a,T]$, it is enough to prove (5.7) without $\pi$. This estimate is elementary. It is enough to split the integral in two pieces, the first one ranging from $-1$ to $-1 + \sqrt{\varepsilon}$, to change variables $v' = \sqrt{\lambda}(1+v)$, and to observe that $\int_0^1 \cosh v\, dv \leq C_1 \int_0^1 \sinh v\, dv$, $\cosh v \leq C_1 \sinh v$ for $v \geq 1$.

An analogous argument shows that there exists a constant $C_1$ such that

$$\int_u^1 dv \cosh[\sqrt{\lambda}(1-v)]\pi(t,v) \leq C_1 \int_u^1 dv \sinh[\sqrt{\lambda}(1-v)]\pi(t,v)$$

for any $0 < \varepsilon \leq 1$, $a \leq t \leq T$, $u$ in $A^\varepsilon$.

By definition of $\pi_\varepsilon$, $\pi_\varepsilon \geq R^D_{\beta_\varepsilon}\pi$. Therefore, by the previous estimate, (5.7) and the explicit form of the kernels $R^N_\varepsilon$, $R^D_\varepsilon$, we easily get that

$$\sup_{0<\varepsilon\leq 1}\sup_{t\in[a,T], u\in A^\varepsilon} \frac{R^N_{\beta_\varepsilon(t)}\pi(t,u)}{\pi_\varepsilon(t,u)} \leq \sup_{0<\varepsilon\leq 1}\sup_{t\in[a,T], u\in A^\varepsilon} \frac{R^N_{\beta_\varepsilon(t)}\pi(t,u)}{R^D_{\beta_\varepsilon(t)}\pi(t,u)} \leq C_0$$

for some finite constant $C_0 = C_0(\pi)$. The same arguments with $\pi$ replaced by $1 - \pi$ yields

$$
\begin{aligned}
& \sup_{0<\varepsilon\leq 1}\sup_{t\in[a,T], u\in A^\varepsilon} \frac{1 - R^N_{\beta_\varepsilon(t)}\pi(t,u)}{1 - \pi_\varepsilon(t,u)} \\
& \qquad \leq \sup_{0<\varepsilon\leq 1}\sup_{t\in[a,T], u\in A^\varepsilon} \frac{R^N_{\beta_\varepsilon(t)}(1-\pi)(t,u)}{R^D_{\beta_\varepsilon(t)}(1-\pi)(t,u)} \leq C_0,
\end{aligned}
$$

which concludes the proof of (5.4). □

Let $a$ be as in the previous proof and note that $\partial_t \pi_\varepsilon$ may be discontinuous at $a$. The left time derivative is equal to $\partial_t \rho$ where $\rho$ is the solution of the hydrodynamic equation (4.2) while the right time derivative vanishes because the derivative of $j$ vanishes at 0 and $\pi$ is constant in time in the interval $[a,b]$. Since $\partial_t \pi_\varepsilon$ vanishes at $a$ for every $\varepsilon > 0$, we could have added this extra assumption in the definition of the set $\mathcal{F}_3 = \mathcal{D}^\circ_{T,\gamma}$, but we do not need it.

Recall the definition of $\mathfrak{t}(\pi)$, given just before the proof of Theorem 5.1.

LEMMA 5.7. *Fix a trajectory $\pi$ in $\mathcal{D}^\circ_{T,\gamma}$. For each $0 \leq t \leq T$, let $H_t$ be the unique solution of the elliptic equation*

(5.8) $$\begin{cases} \partial_t \pi_t = \nabla(D(\pi_t)\nabla\pi_t) - \nabla\{\chi(\pi_t)[(E/2) + \nabla H_t]\}, \\ H_t(\pm 1) = 0. \end{cases}$$

*For $t = \mathfrak{t}(\pi)$, $\partial_t \pi_t$ should be interpreted as the right derivative $\partial_{t+}\pi_t$. Then $H$ vanishes on $[0, \mathfrak{t}(\pi)) \times [-1, 1]$ and $H$ is smooth on $(\mathfrak{t}(\pi), T] \times [-1, 1]$. Moreover,*

$$I_T(\pi|\gamma) = \frac{1}{2}\int_0^T \|H_t\|^2_{1,\chi(\pi_t)}\, dt$$

*and*

$$I_T(\pi|\gamma) = \langle \pi_T, H_T\rangle - \langle \pi_{\mathfrak{t}(\pi)}, H_{\mathfrak{t}(\pi)}\rangle - \int_{\mathfrak{t}(\pi)}^T \langle \pi_t, \partial_t H_t\rangle\, dt$$
$$+ \int_{\mathfrak{t}(\pi)}^T \langle D(\pi_t)\nabla\pi_t, \nabla H_t\rangle\, dt$$
$$- \frac{E}{2}\int_{\mathfrak{t}(\pi)}^T \langle \chi(\pi_t), \nabla H_t\rangle\, dt - \frac{1}{2}\int_{\mathfrak{t}(\pi)}^T \langle \chi(\pi_t), (\nabla H_t)^2\rangle\, dt.$$

PROOF. Fix a trajectory $\pi$ in $\mathcal{D}^\circ_{T,\gamma}$. Since $\pi$ is bounded away from 0 and 1 on $(0, T] \times [-1, 1]$, equation (5.8) is strictly elliptic and, therefore, has a unique solution. Since $\pi$ follows the hydrodynamic equation in the time interval $[0, \mathfrak{t}(\pi)]$, the unique solution $H_t$ is identically equal to 0 on $[0, \mathfrak{t}(\pi)) \times [-1, 1]$. On $[\mathfrak{t}(\pi), T] \times [-1, 1]$ $H$ inherits the smoothness in space from $\pi$. Smoothness in time in this interval follows from the smooth dependence on the force term of solutions of strictly elliptic equations.

It remains to show that the rate function has the explicit forms claimed. This follows from (4.11), or from the next elementary argument. Fix a function $G$ in $C^{1,2}_0(\overline{\Omega_T})$ and recall the definition (4.6) of the functional $\hat{J}_G$. Since $\pi$ solves (5.8), an integration by parts shows that

$$\hat{J}_G(\pi) = \langle\!\langle \chi(\pi)\nabla H, \nabla G\rangle\!\rangle - \tfrac{1}{2}\langle\!\langle \chi(\pi)\nabla G, \nabla G\rangle\!\rangle.$$

Therefore,

$$I_T(\pi|\gamma) = \hat{I}_T(\pi|\gamma)$$
$$= \frac{1}{2}\langle\!\langle \chi(\pi)\nabla H, \nabla H\rangle\!\rangle$$
$$- \frac{1}{2}\inf_{G \in C^{1,2}_0(\overline{\Omega_T})}\{\langle\!\langle \chi(\pi)[\nabla H - \nabla G], [\nabla H - \nabla G]\rangle\!\rangle\}.$$

The last term vanishes because $\chi(\pi)$ is bounded and smooth functions vanishing at the boundary are dense in $L^2(0,T,H_0^1(\Omega))$.

To prove the second identity, multiply (5.8) by $H$, integrate in space and in time in the interval $[\mathfrak{t}(\pi), T]$, and integrate by parts to get that

$$\langle \pi_T, H_T\rangle - \langle \pi_{\mathfrak{t}(\pi)}, H_{\mathfrak{t}(\pi)}\rangle - \int_{\mathfrak{t}(\pi)}^T \langle \pi_t, \partial_t H_t\rangle\, dt + \int_{\mathfrak{t}(\pi)}^T \langle D(\pi_t)\nabla\pi_t, \nabla H_t\rangle\, dt$$
$$= \frac{E}{2}\int_{\mathfrak{t}(\pi)}^T \langle \chi(\pi_t), \nabla H_t\rangle\, dt + \int_{\mathfrak{t}(\pi)}^T \langle \chi(\pi_t), (\nabla H_t)^2\rangle\, dt.$$

Since $H$ vanishes in the time interval $[0,\mathfrak{t}(\pi)]$, we may replace in the last integral $\mathfrak{t}(\pi)$ by 0. The second formula for the rate function $I_T(\pi|\gamma)$ follows from this observation, the previous identity and the first formula for $I_T(\pi|\gamma)$. □

**Acknowledgments.** The authors thank the referees for their careful reading of the manuscript and their comments which helped to improve presentation of the article.

## REFERENCES

[1] BENOIS, O. (1996). Large deviations for the occupation times of independent particle systems. *Ann. Appl. Probab.* **6** 269–296. MR1389840
[2] BENOIS, O., KIPNIS, C. and LANDIM, C. (1995). Large deviations from the hydrodynamical limit of mean zero asymmetric zero range processes. *Stochastic Process. Appl.* **55** 65–89. MR1312149
[3] BERTINI, L., DE SOLE, A., GABRIELLI, D., JONA-LASINIO, G. and LANDIM, C. (2003). Large deviations for the boundary driven symmetric simple exclusion process. *Math. Phys. Anal. Geom.* **6** 231–267. MR1997915
[4] BERTINI, L., DE SOLE, A., GABRIELLI, D., JONA-LASINIO, G. and LANDIM, C. (2006). Large deviation approach to nonequilibrium processes in stochastic lattice gases. *Bull. Braz. Math. Soc. (N.S.)* **37** 611–643. MR2284891
[5] BERTINI, L., GABRIELLI, D. and LANDIM, C. (2009). Strong asymmetric limit of the quasi-potential of the boundary driven weakly asymmetric exclusion process. *Comm. Math. Phys.* **289** 311–334.
[6] COMETS, F. (1986). Grandes déviations pour des champs de Gibbs sur $Z^d$. *C. R. Acad. Sci. Paris Sér. I Math.* **303** 511–513. MR865873
[7] DE MASI, A., PRESUTTI, E. and SCACCIATELLI, E. (1989). The weakly asymmetric simple exclusion process. *Ann. Inst. H. Poincaré Probab. Statist.* **25** 1–38. MR995290
[8] DERRIDA, B. (2007). Nonequilibrium steady states: Fluctuations and large deviations of the density and of the current. *J. Stat. Mech. Theory Exp.* **7** 45. MR2335699
[9] DERRIDA, B., LEBOWITZ, J. L. and SPEER, E. R. (2002). Large deviation of the density profile in the steady state of the open symmetric simple exclusion process. *J. Stat. Phys.* **107** 599–634. MR1898851
[10] DERRIDA, B., LEBOWITZ, J. L. and SPEER, E. R. (2003). Exact large deviation functional of a stationary open driven diffusive system: The asymmetric exclusion process. *J. Stat. Phys.* **110** 775–810. MR1964689

[11] DEUSCHEL, J.-D. and STROOCK, D. W. (1989). *Large Deviations. Pure and Applied Mathematics* **137**. Academic Press, Boston, MA. MR997938

[12] DONSKER, M. D. and VARADHAN, S. R. S. (1989). Large deviations from a hydrodynamic scaling limit. *Comm. Pure Appl. Math.* **42** 243–270. MR982350

[13] ENAUD, C. and DERRIDA, B. (2004). Large deviation functional of the weakly asymmetric exclusion process. *J. Stat. Phys.* **114** 537–562. MR2035624

[14] EYINK, G., LEBOWITZ, J. L. and SPOHN, H. (1990). Hydrodynamics of stationary nonequilibrium states for some stochastic lattice gas models. *Comm. Math. Phys.* **132** 253–283. MR1069212

[15] EYINK, G., LEBOWITZ, J. L. and SPOHN, H. (1991). Lattice gas models in contact with stochastic reservoirs: Local equilibrium and relaxation to the steady state. *Comm. Math. Phys.* **140** 119–131. MR1124262

[16] GÄRTNER, J. (1988). Convergence towards Burgers' equation and propagation of chaos for weakly asymmetric exclusion processes. *Stochastic Process. Appl.* **27** 233–260. MR931030

[17] KIPNIS, C. and LANDIM, C. (1999). *Scaling Limits of Interacting Particle Systems*. Springer, Berlin. MR1707314

[18] KIPNIS, C., LANDIM, C. and OLLA, S. (1995). Macroscopic properties of a stationary nonequilibrium distribution for a nongradient interacting particle system. *Ann. Inst. H. Poincaré Probab. Statist.* **31** 191–221. MR1340037

[19] KIPNIS, C., MARCHIORO, C. and PRESUTTI, E. (1982). Heat flow in an exactly solvable model. *J. Stat. Phys.* **27** 65–74. MR656869

[20] KIPNIS, C., OLLA, S. and VARADHAN, S. R. S. (1989). Hydrodynamics and large deviations for simple exclusion processes. *Comm. Pure Appl. Math.* **42** 115–137. MR978701

[21] LANDIM, C., OLLA, S. and VOLCHAN, S. B. (1998). Driven tracer particle in one-dimensional symmetric simple exclusion. *Comm. Math. Phys.* **192** 287–307. MR1617558

[22] LANFORD, O. E. (1973). *Entropy and Equilibrium States is Classical Statistical Mechanics. Lecture Notes in Physics* **20**. Springer, Berlin.

[23] MOURRAGUI, M. and ORLANDI, E. (2007). Large deviations from a macroscopic scaling limit for particle systems with Kac interaction and random potential. *Ann. Inst. H. Poincaré, Probab. Statist.* **43** 677–715.

[24] OLLA, S. (1988). Large deviations for Gibbs random fields. *Probab. Theory Related Fields* **77** 343–357. MR931502

[25] PROTTER, M. H. and WEINBERGER, H. F. (1984). *Maximum Principles in Differential Equations*. Springer, New York. MR762825

[26] QUASTEL, J. (1995). Large deviations from a hydrodynamic scaling limit for a nongradient system. *Ann. Probab.* **23** 724–742. MR1334168

[27] QUASTEL, J., REZAKHANLOU, F. and VARADHAN, S. R. S. (1999). Large deviations for the symmetric simple exclusion process in dimensions $d \geq 3$. *Probab. Theory Related Fields* **113** 1–84. MR1670733

[28] SIMON, J. (1987). Compact sets in the space $L^p(0,T;B)$. *Ann. Mat. Pura Appl. (4)* **146** 65–96. MR916688

[29] SPOHN, H. and YAU, H.-T. (1995). Bulk diffusivity of lattice gases close to criticality. *J. Stat. Phys.* **79** 231–241.

[30] VARADHAN, S. R. S. (1984). *Large Deviations and Applications*. SIAM, Philadelphia. MR758258

[31] VARADHAN, S. R. S. and YAU, H.-T. (1997). Diffusive limit of lattice gas with mixing conditions. *Asian J. Math.* **1** 623–678. MR1621569

[32] Zeidler, E. (1990). *Nonlinear Functional Analysis and Its Applications. II/A*. Springer, Berlin. MR1033497

L. Bertini  
Dipartimento di Matematica  
Università di Roma "La Sapienza"  
P. le Aldo Moro 2, 00185 Roma  
Italy  
E-mail: bertini@mat.uniroma1.it

C. Landim  
IMPA  
Estrada Dona Castorina 110  
J. Botanico, 22460 Rio de Janeiro  
Brazil  
E-mail: landim@impa.br

C. Landim  
M. Mourragui  
CNRS UMR 6085  
Université de Rouen  
Avenue de l'Université, BP.12  
Technopôle du Madrillet  
F76801 Saint-Étienne-du-Rouvray  
France  
E-mail: Mustapha.Mourragui@univ-rouen.fr